\documentclass[11pt]{amsart}

\usepackage{amssymb}

\usepackage{upref}

\newcommand{\conservepaper}{
 \hoffset=-0.75in
 \setlength{\textwidth}{6.5in}
 \voffset=-0.5in
 \setlength{\textheight}{10.0in}
 \setlength{\textheight}{9.0in} 
 }
 \conservepaper

\swapnumbers 

\numberwithin{equation}{section}
\newtheorem{thm}[equation]{Theorem}

\newtheorem{lem}[equation]{Lemma}
\newtheorem{cor}[equation]{Corollary}
\newtheorem{prop}[equation]{Proposition}

\theoremstyle{definition}

\newtheorem{defn}[equation]{Definition}

\newtheorem{eg}[equation]{Example}

\newtheorem{notation}[equation]{Notation}

\theoremstyle{remark}

\newtheorem{rem}[equation]{Remark}        
         
\newtheorem{claim}[equation]{Claim}
\newtheorem{ack}[equation]{Acknowledgment}  

\newcommand{\pref}[1]{{\upshape(}\ref{#1}{\upshape)}}
\renewcommand{\see}[1]{{\upshape(}see~\ref{#1}{\upshape)}}
\newcommand{\cf}[1]{(c.f.~\ref{#1})}
\newcommand{\fullref}[2]{\ref{#1}\pref{#1#2}}

\newcommand{\Lie}[1]{\mathfrak{#1}}
\newcommand{\diag}{\operatorname{diag}}
\newcommand{\transpose}{T}
\newcommand{\SL}{\operatorname{SL}}
\newcommand{\GL}{\operatorname{GL}}
\newcommand{\SO}{\operatorname{SO}}
\newcommand{\so}{\operatorname{\Lie{so}}}

\newcommand{\Id}{\operatorname{Id}}

\newcommand{\Rad}{\operatorname{Rad}}
\newcommand{\real}{\mathord{\mathbb{R}}}

\newcommand{\integer}{\mathord{\mathbb{Z}}}

\newcommand{\semiprod}{\ltimes}
\newcommand{\closure}[1]{\overline{#1}}

\newcommand{\Rrank}{\mathop{\real\text{\upshape-rank}}}
\newcommand{\chamber}{\mathord{\mathcal{C}}}

\newcommand{\bdry}{\partial}

\newcommand{\muH}[2]{\bigl[ #1, #2 \bigr]}

\newcommand{\bigset}[2]{\left\{\, #1 
 \mathrel{\left| \vphantom {\left\{ #1 \mid #2 \right\} } \right.}
 #2 \,\right\} }

\newcommand{\h}[6]{\begin{pmatrix}
 {#1}&{#4}&{#6}\\
 0&{#2}&{#5}\\
 0&0&{#3}\\
 \end{pmatrix}
 }

 \newcounter{case}
 \newenvironment{case}[1][\unskip]{\refstepcounter{case}\em
 \medskip \noindent Case \thecase\ #1.\ }{\unskip\upshape}
 \renewcommand{\thecase}{\arabic{case}}

 \newcounter{subcase}
 \newenvironment{subcase}[1][\unskip]{\refstepcounter{subcase}\em
 \medskip \noindent Subcase \thesubcase\ #1.\ }{\unskip\upshape}
\numberwithin{subcase}{case}

\renewenvironment{claim}[1][\unskip]{
\em
 \medskip \noindent Claim
 #1.\ }{\unskip\upshape}

\begin{document}

\title{Cartan-decomposition subgroups of $\SO(2,n)$}

\author{Hee Oh}
 \address{Department of Mathematics, Oklahoma State University,
Stillwater, OK 74078}
 \curraddr{Institute of Mathematics, The Hebrew University,
 Jerusalem 91904, Israel}
 \email{heeoh@math.huji.ac.il}

\author{Dave Witte}
 \address{Department of Mathematics, Oklahoma State University,
Stillwater, OK 74078}
 \email{dwitte@math.okstate.edu,
  http://www.math.okstate.edu/\char'176dwitte} 


\date{February 26, 1999 \bf (Corrected version)} 

\begin{abstract}
 For $G = \SL(3,\real)$ and $G = \SO(2,n)$, we give explicit, practical
conditions that determine whether or not a closed, connected subgroup~$H$
of~$G$ has the property that there exists a compact subset~$C$ of~$G$
with $CHC = G$. To do this, we fix a Cartan decomposition $G = K A^+ K$
of~$G$, and then carry out an approximate calculation of $(KHK) \cap A^+$
for each closed, connected subgroup~$H$ of~$G$.
 \end{abstract}

\maketitle

\section{Introduction}

\begin{notation}
 Throughout this paper, $G$ is a Zariski-connected, almost simple,
linear, real Lie group. (``Almost simple" means that every proper normal
subgroup of~$G$ either is finite or has finite index.) In almost all of
the main results, $G$ is assumed to be either $\SL(3,\real)$ or
$\SO(2,n)$ (with $n \ge 3$). There would be no essential loss of
generality if one were to require $G$ to be connected, instead of only
Zariski connected \see{Gconn}. However, $\SO(2,n)$ is not connected (it
has two components) and the authors prefer to state results for
$\SO(2,n)$, instead of for the identity component of $\SO(2,n)$.

Fix an Iwasawa decomposition  $G = KAN$ and a corresponding Cartan
decomposition $G = K A^+ K$, where $A^+$ is the (closed) positive Weyl
chamber of~$A$ in which the roots occurring in the Lie algebra of~$N$ are
positive. Thus, $K$ is a maximal compact subgroup, $A$ is the identity
component of a maximal split torus, and $N$ is a maximal unipotent
subgroup.
 \end{notation}

The terminology introduced in the following definition is new, but the
underlying concept is well known (see, for example,
Proposition~\ref{Calabi-Markus} and Theorem~\ref{CDSvsmu} below).

\begin{defn}
 Let $H$ be a closed subgroup of~$G$. We say that $H$ is a
\emph{Cartan-decomposition subgroup} of~$G$ if 
 \begin{itemize}
 \item $H$ is connected, and
 \item there is a compact
subset~$C$ of~$G$, such that $CHC = G$. 
 \end{itemize}
 (Note that $C$ is only assumed to
be a sub\emph{set} of~$G$; it need not be a sub\emph{group}.)
 \end{defn}

\begin{eg}
 The Cartan decomposition $G = KAK$ shows that the maximal split torus~$A$
is a Cartan-decomposition subgroup of~$G$.

It is known that $G = KNK$ \cite[Thm.~5.1]{Kostant}, so the maximal
unipotent subgroup~$N$ is also a Cartan-decomposition subgroup.

 If $G$ is compact (that is, if $\Rrank G = 0$), then all subgroups
of~$G$ are Cartan-decomposition subgroups. 
 On the other hand, if $G$ is noncompact, then not all subgroups are
Cartan-decomposition subgroups, because it is obvious that every
Cartan-decomposition subgroup of~$G$ is noncompact. It is somewhat less
obvious that if $H$ is a Cartan-decomposition subgroup of~$G$, then $\dim
H \ge \Rrank G$ \see{dim>Rrank}.
 \end{eg}

Our interest in  Cartan-decomposition subgroups is largely motivated by
the following basic observation that, to construct nicely behaved actions
on homogeneous spaces, one must find subgroups that are \emph{not}
Cartan-decomposition subgroups. (See \cite[\S3]{Kobayashi-survey} for
some historical background on this result.)

\begin{prop}[{Calabi-Markus phenomenon, cf.~\cite[pf.~of
Thm.~A.1.2]{Kulkarni}}] \label{Calabi-Markus}
 If $H$ is a Cartan-decomposition subgroup of~$G$, then no closed,
noncompact subgroup of~$G$ acts properly on $G/H$.
 \end{prop}

Our goal is to determine which closed, connected subgroups of~$G$ are
Cartan-decomposition subgroups, and which are not. Our main tool is the
Cartan projection.

\begin{defn}[Cartan projection]
 For each element~$g$ of~$G$, the Cartan decomposition $G = K A^+ K$
implies that there is an element~$a$ of~$A^+$ with $g \in K a K$. In fact,
the element~$a$ is unique, so there is a well-defined function $\mu
\colon G \to A^+$ given by $g \in K \, \mu(g) \, K$. The function $\mu$ is
continuous and proper (that is, the inverse image of any compact set
is compact). Some properties of the Cartan projection are discussed
in~\cite{Benoist} and~\cite{Kobayashi-survey}.
 \end{defn}

We have $\mu(H) = A^+$ if and only if $KHK = G$. This immediately implies
that if $\mu(H) = A^+$, then $H$ is a Cartan-decomposition subgroup.
Y.~Benoist and T.~Kobayashi proved the deeper statement that, in the
general case, $H$ is a Cartan-decomposition subgroup if and only if
$\mu(H)$ comes within a bounded distance of every point in~$A^+$.

\begin{thm}[Benoist {\cite[Prop.~5.1]{Benoist}, Kobayashi
\cite[Thm.~1.1]{Kobayashi-criterion}}] \label{CDSvsmu}
 \ 
 A closed, connected subgroup~$H$ of~$G$ is a Cartan-decomposition
subgroup if and only if there is a compact subset~$C$ of~$A$, such that
$\mu(H) C \supset A^+$.
 \end{thm}

 We noted above that every subgroup is a Cartan-decomposition subgroup
if $\Rrank G = 0$. Therefore, the characterization of
Cartan-decomposition subgroups of~$G$ is trivial if $\Rrank G = 0$. The
following simple proposition shows that the characterization is again
very easy if $\Rrank G = 1$.

\begin{prop}[{cf.~\cite[Lem.~3.2]{Kobayashi-isotropy}}]
\label{Rrank1-CDS}
 Assume that $\Rrank G = 1$. A closed, connected subgroup~$H$ of~$G$ is a
Cartan-decomposition subgroup if and only if $H$ is noncompact.
 \end{prop}

\begin{proof}
 ($\Leftarrow$) We have $\mu(e) = e$, and, because $\mu$ is a proper map,
we have $\mu(h) \to \infty$ as $h \to \infty$ in~$H$. Because $\Rrank G =
1$, we know that $A^+$ is homeomorphic to the half-line $[0,\infty)$
(with the point~$e$ in~$A^+$ corresponding to the endpoint~$0$ of the
half-line), so, by continuity, it must be the case that $\mu(H) = A^+$.
Therefore $KHK = G$, so $H$ is a Cartan-decomposition subgroup.
 \end{proof}

It seems to be much more difficult to characterize the
Cartan-decomposition subgroups when $\Rrank G = 2$, so these are the first
interesting cases. In this paper, we study two examples in detail.
Namely, we describe all the Cartan-decomposition subgroups of
$\SL(3,\real)$ and of $\SO(2,n)$. We also explicitly describe the closed,
connected subgroups that are \emph{not} Cartan-decomposition subgroups,
and approximately calculate the image of each of these subgroups under the
Cartan projection.

Obviously, any connected, closed subgroup that contains a
Cartan-decomposition subgroup is itself a Cartan-decomposition subgroup.
Therefore the minimal Cartan-decomposition subgroups are the most
interesting ones. As a simple example of our results, we state the
following theorem.

\begin{thm} \label{SL3-minl}
 Assume that $G = \SL(3,\real)$.
 Up to conjugation by automorphisms of~$G$, the only minimal
Cartan-decomposition subgroups of~$G$ are:
 $$ A,
  \qquad
 \bigset{\h111rrs}{r,s \in \real},
  \qquad
 \bigset{\h{e^t}{e^t}{e^{-2t}}{te^t}rs}{r,s,t \in \real},
 $$
 and subgroups of the form
 \begin{equation} \label{SL3-rootsemi}
 \bigset{\h{e^{pt}}{e^{qt}}{e^{-(p+q)t}}r00}{r,t \in \real}, 
 \end{equation}
 where $p$~and~$q$ are fixed real numbers with
$\max\{p,q\} = 1$ and $\min\{p,q\} \ge -1/2$, or of the form
 \begin{equation} \label{SL3-rot}
  \bigset{
 \begin{pmatrix}
   e^{t} \cos pt   & e^{t} \sin pt  & s \\
 - e^{t} \sin pt   & e^{t} \cos pt  & r \\
   0               &    0           & e^{-2t} \\
 \end{pmatrix}
 }{r,s,t \in \real}, 
 \end{equation}
 where $p$ is a fixed nonzero real number.
 \end{thm}

Note that $AN$ contains uncountably many nonconjugate
minimal Cartan-decomposition subgroups of~$G$, because the minimum of
the two parameters $p$~and~$q$ in~\pref{SL3-rootsemi} can be varied
continuously. However, up to conjugacy under $\operatorname{Aut} G$,
there is only one minimal Cartan-decomposition subgroup contained in~$A$
(namely, $A$~itself), and only one contained in~$N$.

\begin{cor}
 Let $H$ be a closed, connected subgroup of $G = \SL(3,\real)$, and let
$K_H$ be a maximal compact subgroup of~$H$. If $\dim H - \dim K_H \ge 3$,
then $H$ is a Cartan-decomposition subgroup of~$G$.
 \end{cor}

For an explicit description of the Cartan-decomposition subgroups of
$\SL(3,\real)$ (not just up to conjugacy), see Theorem~\ref{SL3-CDS}. The
subgroups of $\SL(3,\real)$ that are not Cartan-decomposition subgroups
are described in Corollary~\ref{SL3notCDS}, and their images under the
Cartan projection are described in Proposition~\ref{SL3-Cartanproj}.
These results are stated only for subgroups of~$AN$, because the general
case reduces to this (see Remark~\ref{CDS-onlyAN}).

Theorem~\ref{SO2n-inN-minl} is a sample of our results on
Cartan-decomposition subgroups of $\SO(2,n)$. Note that, for simplicity,
we restrict here to subgroups of~$N$. 

\begin{notation} \label{SO2n-defn}
 We realize $\SO(2,n)$ as isometries of the indefinite form $\langle v
\mid v \rangle = v_1 v_{n+2} + v_2 v_{n+1} + \sum_{i=3}^n v_i^2$
on~$\real^{n+2}$ (for $v = (v_1,v_2,\ldots,v_{n+2}) \in \real^{n+2}$).
 The virtue of this particular realization is that we may choose $A$ to
consist of the diagonal matrices in $\SO(2,n)$ (with nonnegative
entries) and $N$ to consist of the upper-triangular matrices in
$\SO(2,n)$ with only $1$'s on the diagonal.
 Thus, the Lie algebra of $AN$ is
 \begin{equation}
 \label{SO2n-AN}
 \Lie a + \Lie n =
 \bigset{
  \begin{pmatrix}
 t_1 & \phi & x  & \eta &0\\
 & t_2 & y &0 & -\eta \\
  &   & 0 & -y^T & -x^T \\
 && &-t_2& -\phi \\
 &&& &-t_1 \\
 \end{pmatrix}
 }{ {t_1,t_2,\phi,\eta \in \real, \atop x,y \in \real^{n-2}}}
 .
 \end{equation}
 Note that  the first two rows of any element of $\Lie a + \Lie n$ are
sufficient to determine the entire matrix.
 \end{notation}

Whenever $m < n$, there is an obvious embedding of $\SO(2,m)$ in
$\SO(2,n)$, induced by an inclusion $\real^{m+2} \hookrightarrow
\real^{n+2}$, so, abusing notation, we speak of $\SO(2,m)$ as a
subgroup of $\SO(2,n)$.

 \begin{thm} \label{SO2n-inN-minl}
 Assume that $G = \SO(2,n)$. If $n \ge 5$, then there are exactly~6
non-conjugate minimal Cartan-decomposition subgroups of~$G$ contained
in~$N$. Each such subgroup~$H$ is conjugate to a subgroup of $\SO(2,5)$
and, as a subalgebra of $\so(2,5)$, the Lie algebra of~$H$ is conjugate
to one of the following:
 \begin{enumerate}
 \item \label{SO2n-inN-minl-a,a+2b}
 $
 \bigset{
 \begin{pmatrix}
 0&\phi &0               &0&0 &\eta&0 \\
  &0    &\epsilon_1 \phi &0&0 &0&-\eta \\
  &&& \cdots \\
 \end{pmatrix}
  }{ \phi,\eta \in \real}$,
 where $\epsilon_1 \in \{0,1\}$
 \item \label{SO2n-inN-minl-a,a+b}
 $
 \bigset{
 \begin{pmatrix}
 0&\phi &x&0&0                &0 &0 \\
  &0    &0& \epsilon_2 \phi&0 &0&0 \\
  &&& \cdots \\
 \end{pmatrix}
 }{ \phi, x \in \real}$,
  where $\epsilon_2 \in \{0,1\}$
 \item \label{SO2n-inN-minl-a+b,b}
 $
 \bigset{
 \begin{pmatrix}
 0&0 &x&0&\epsilon_3 y &0&0 \\
  &0 &0&y&0          &0&0 \\
  &&& \cdots \\
 \end{pmatrix}
 }{ x,y \in \real}
 $,
 where $\epsilon_3 \in \{0,1\}$.
 \end{enumerate}

 There are 5 non-conjugate minimal Cartan-decomposition
subgroups of~$\SO(2,4)$ contained in~$N$. The Lie algebra of any such
subgroup~$H$ is conjugate either to one of the two
subalgebras of type~\pref{SO2n-inN-minl-a,a+2b}, 
to one of the two subalgebras of type~\pref{SO2n-inN-minl-a,a+b}, or to
the subalgebra of type~\pref{SO2n-inN-minl-a+b,b} with $\epsilon_3 = 0$.
{\upshape(}These are the five of the above-listed subalgebras that are
contained in $\so(2,4)$, namely, the five whose 5th column is
all~$0$'s.{\upshape)}

 There are 3 non-conjugate minimal Cartan-decomposition
subgroups of~$\SO(2,3)$ contained in~$N$. The Lie algebra of any such
subgroup~$H$ is conjugate  either to one of the two
subalgebras of type~\pref{SO2n-inN-minl-a,a+2b},  or to the subalgebra of
type~\pref{SO2n-inN-minl-a,a+b} with $\epsilon_2 = 0$. {\upshape(}These
are the three of the above-listed subalgebras that are contained in
$\so(2,3)$, namely, the three whose 4th and~5th columns are
all~$0$'s.{\upshape)}
 \end{thm}

The detailed study of Cartan-decomposition subgroups of $\SO(2,n)$ is
rather complicated, so we break it up into three parts: subgroups of~$N$
(Theorem~\ref{SO2n-HinN-CDS}), subgroups not in~$N$ that can be written
as a semidirect product $T \semiprod U$ with $T \subset A$ and $U \subset
N$ (Theorem~\ref{SO2n-semiprod}), and subgroups that cannot be written as
such a semidirect product (Theorem~\ref{SO2n-notsemi-CDS}). We also
describe the subgroups of $\SO(2,n)$ that are not Cartan-decomposition
subgroups (see Theorem~\ref{HinN-notCDS} and
Corollaries~\ref{SO2n-semi-notCDS} and~\ref{SO2n-notsemi-notCDS}), and
approximately calculate their Cartan projections (see
Proposition~\ref{HinN-Cartanproj} and Corollaries~\ref{SO2n-semi-notCDS}
and~\ref{SO2n-notsemi-notCDS}).

If $H$ is a Cartan-decomposition subgroup of~$G$, and $G/H$ is not
compact, then the Calabi-Markus phenomenon~\ref{Calabi-Markus} implies
that $G/H$ does not have a compact Clifford-Klein form.  (That is, there
does not exist a discrete subgroup~$\Gamma$ of~$G$ that acts properly on
$G/H$, such that the quotient space $\Gamma \backslash G/H$ is compact.)
Thus, our work on Cartan-decomposition subgroups is a first step toward
understanding which homogeneous spaces of~$G$ have a compact
Clifford-Klein form. Building on this, a sequel~\cite{OhWitte-CK}
determines exactly which homogeneous spaces of $\SO(2,n)$ have a compact
Clifford-Klein form in the case where $n$ is even (and assuming that the
isotropy group~$H$ is connected), but the results are not quite complete
when $n$ is odd. The work leads to new examples of compact
Clifford-Klein forms of $\SO(2,n)$, when $n$ is even.

The paper is organized as follows.
 Section~\ref{prelim} collects some
known results on Lie groups and Zariski closures.
 Section~\ref{CDS-general} presents some general results on
Cartan-decomposition subgroups.
 Section~\ref{SL3R} states and proves
our results on Cartan-decomposition subgroups of $\SL(3,\real)$.
 Section~\ref{SO2n-HinN-sect} contains our results
on Cartan-decomposition subgroups of $\SO(2,n)$ that are contained
in~$N$, and
 Section~\ref{SO2n-notinN-sect} is devoted to the subgroups of $\SO(2,n)$
that are not contained in~$N$.

\begin{ack} This research was partially supported by grants from the
National Science Foundation (DMS-9623256 and DMS-9801136). 
 Many of the results of this paper were finalized during a visit of the
authors to the University of Bielefeld. We would like to  thank the
German-Israeli Foundation for Research and Development for financial
support that made the visit possible, and the mathematics faculty of the
University of Bielefeld for their hospitality that helped make the visit
so productive.  We are also grateful to Kari Vilonen for explaining that
the theory of o-minimal structures would be useful to us
(see~\ref{dim>Rrank} and~\ref{CDS<>seq}). D.~Witte would like to thank
the Tata Institute for Fundamental Research for providing a congenial
environment to carry out final revisions on the manuscript.
 \end{ack}

\section{Preliminaries on Lie groups and Zariski closures} \label{prelim}

Most of the results in this section are well known, and none are new. The
reader is encouraged to skip over this section, and refer back when
necessary.

We assume familiarity with the basic theory of Lie groups and Lie
algebras (as in, for example, \cite{Hochschild-Lie}). At some points, we
also assume some familiarity with the structure of algebraic groups
over~$\real$, in the spirit of \cite[\S P.2, pp.~7--11]{Raghunathan}.

\begin{notation}
 We use German letters $\Lie g$, $\Lie h$, $\Lie a$, $\Lie n$, $\Lie u$,
$\Lie t$ to denote the Lie algebras of Lie groups $G$, $H$, $A$, $N$, $U$,
$T$, etc.

 For a linear functional~$\omega$ on~$\Lie a$, we use $\Lie u_\omega$ to
denote the corresponding weight space of the Lie algebra~$\Lie g$ (so
$\Lie u_\omega = 0$ unless $\omega$ is either~$0$ or a real
root of~$G$), and $U_\omega$ to denote $\exp \Lie u_\omega$. Note that
$U_\omega U_{2\omega}$ is a subgroup of~$G$ (but $U_\omega$ may not be a
subgroup if both $\omega$ and~$2\omega$ are real roots of~$G$).
 \end{notation}

The idea of the following definition is to require that the choice of the
Cartan subgroup~$A$ be compatible with a particular subgroup~$H$. It is
not a severe restriction on~$H$, because Lemma~\ref{conj-to-compatible}
shows that it can always be satisfied by replacing~$H$ with a conjugate.
Then, under the assumption that $\Rrank G = 2$, Lemma~\ref{not-semi}
shows that $H$ has a fairly simple description in terms of (1)~an
appropriate subgroup of~$A$, (2)~the intersection of~$H$ with~$N$, and,
perhaps, (3)~a~homomorphism into a root subgroup of~$G$. (If $\Rrank G >
2$, then the description of a typical~$H$ would require several root
subgroups.)

\begin{defn} \label{compatible}
 Let us say that a subgroup~$H$ of $AN$ is \emph{compatible} with~$A$ if
$H \subset T U C_N(T)$, where $T = A \cap (HN)$, $U = H \cap N$, and
$C_N(T)$ denotes the centralizer of~$T$ in~$N$. 
 \end{defn}

\begin{lem} \label{conj-to-compatible}
 If $H$ is a closed, connected subgroup of~$AN$, then $H$ is conjugate,
via an element of~$N$, to a subgroup that is compatible with~$A$.
 \end{lem}

\begin{proof}
 Let $\closure{H}$ be the identity component of the Zariski closure
of~$H$, and write $\closure{H} = \closure{T} \semiprod \closure{U}$,
where $\closure{U}$ is a subgroup of~$N$ and $\closure{T}$ is conjugate,
via an element of~$N$, to a subgroup of~$A$. Replacing $H$ by a
conjugate, we may assume that $\closure{T} \subset A$. Let $U = H \cap
N$. Then $[\closure{H},\closure{H}] \subset H \cap N$
(cf.~\cite[Cor.~3.8.4, p.~207]{Varadarajan} and \cite[Lem~3.24]{Witte}),
so we have $[\closure{T},\closure{\Lie u}] \subset \Lie u$. Because
$\closure{T} \subset A$ and $\closure{\Lie u}$ is
$(\operatorname{Ad}_G(T))$-invariant, the adjoint action
of~$\closure{T}$ on~$\closure{\Lie u}$ is completely reducible, so this
implies that there is a subspace $\Lie c$ of~$\closure{\Lie u}$, such
that $[\closure{T}, \Lie c] = 0$ and $\Lie u + \Lie c = \closure{\Lie
u}$. Therefore, $U \, C_N(\closure{T}) = \closure{U}$, so
 $\closure{H} = \closure{T} U C_N(\closure{T})$.
 
 Let $\pi \colon AN \to A$ be the projection with kernel~$N$,
and let $T = \pi(H)$. Then $T = \pi(H) \subset \pi(\closure{H}) =
\closure{T}$, so $C_N(T) \supset C_N(\closure{T})$. For any $h \in H$,
there exist $t \in \closure{T}$, $u \in U$ and $c \in C_N(\closure{T})$,
such that $h = tuc$. Because $uc \in N$, we must have $t = \pi(h) \in T$
and, because $C_N(T) \supset C_N(\closure{T})$, we have $c \in C_N(T)$.
Therefore, $h \in T U C_N(T)$. We conclude that $H \subset T U C_N(T)$,
so $H$ is compatible with~$A$.
 \end{proof}

\begin{lem} \label{not-semi}
 Assume that $\Rrank G = 2$.
 Let $H$ be a closed, connected subgroup of~$AN$, and assume that $H$ is
compatible with~$A$. Then either
 \begin{enumerate}
 \item $H = (H \cap A) \semiprod (H \cap N)$; or
 \item \label{not-semi-not}
 there is a positive root~$\omega$, a nontrivial group homomorphism
$\psi\colon \ker \omega \to U_\omega U_{2\omega}$, and a closed, connected
subgroup~$U$ of~$N$, such that
 \begin{enumerate}
 \item \label{not-semi-codim1}
 $H = \{\, a \psi(a) \mid a \in \ker \omega \,\} U$;
 \item \label{not-semi-disjoint}
 $U \cap \psi( \ker\omega) = e$; and
 \item \label{not-semi-normal}
 $U$ is normalized by both $\ker\omega$ and~$\psi( \ker\omega)$.
 \end{enumerate}
 \end{enumerate}
 \end{lem}

\begin{proof} 
 Because $H$ is compatible with~$A$, we have $H \subset T U C_N(T)$,
where $T = A \cap (HN)$ and $U = H \cap N$. We may assume that $H \neq TU$,
for otherwise we have $H = (H \cap A) \semiprod (H \cap N)$. Therefore
$C_N(T) \neq e$. Because $\Lie n$ is a sum of root spaces, this implies
that there is a positive root~$\omega$, such that $T \subset \ker\omega$.
Because $\Rrank G = 2$, we must have $T =\ker\omega$, for otherwise we
would have $T = e$, so $H = U = TU$. Therefore, $C_N(T) = U_\omega
U_{2\omega}$. 

Because $U \subset H \subset T U C_N(T)$, we have $H = U \bigl[H \cap
\bigl(T C_N(T) \bigr) \bigr]$, so there is a nontrivial one-parameter
subgroup $\{x^t\}$ in $H \cap \bigl(T C_N(T) \bigr)$ that is not
contained in~$U$. Because $T$ centralizes $C_N(T)$, we may write $x^t =
a^t u^t$ where $\{a^t\}$ is a one-parameter subgroup of~$T$ and $\{u^t\}$
is a one-parameter subgroup of $C_N(T)$. Furthermore, this Jordan
decomposition is unique, because $T \cap C_N(T) = e$. Replacing $H$ by a
conjugate subgroup, we may assume that $a^t \in A$. Define $\psi \colon
\ker\omega \to U_\omega U_{2\omega}$ by $\psi(a^t) = u^t$ for all $t \in
\real$. 

\pref{not-semi-codim1} For all $t \in
\real$, we have $a^t \psi(a^t) = a^t u^t = x^t \in H$, which establishes
one inclusion of~\pref{not-semi-codim1}. The other will follow if we show
that $\dim H - \dim U = 1$, so suppose $\dim H - \dim U \ge 2$. Then
Lemma~\ref{HN=AN} implies that $A \subset H$, so it follows from
Lemma~\ref{normedbyA} that $H = A \semiprod (H \cap N)$, contradicting
our assumption that $H \neq TU$.

\pref{not-semi-disjoint} Suppose $U \cap \psi(\ker\omega) \neq e$.
Because the exponential map from~$\Lie n$ to~$N$ is bijective, we know
that the intersection is connected. Because $\dim(\ker\omega) = 1$, this
implies that $\psi(\ker\omega) \subset U$. Therefore $a^t = x^t u^{-t}
\in HU = H$, so $T \subset H$. This contradicts the fact that $H \neq
TU$.

 \pref{not-semi-normal} Because $x^t \in H$, we know that the Jordan
components $a^t$ and~$u^t$ of~$x^t$ belong to the Zariski closure of~$H$
\cite[Thm.~15.3, p.~99]{Humphreys-AlgicGrp}. Therefore, both of $a^t$ and~$u^t$
normalize~$H$ \see{Zar-norm}. Being in $AN$, they also normalize~$N$.
Therefore, they normalize $H \cap N = U$.
 \end{proof}

\begin{lem}[{cf.~\cite[pf.~of Thm.~20.2(d),
pp.~108--109]{Humphreys-LieAlg}}] \label{normedbyA}
 If $H$ is a closed connected subgroup of~$AN$ that is normalized
by~$A$, and $\omega$ is a weight of the adjoint representation of~$A$
on~$\Lie a + \Lie n$, then $\pi_\omega(\Lie h) \subset \Lie h$, where
$\pi_\omega \colon \Lie a + \Lie n \to \Lie u_\omega$ is the
$A$-equivariant projection.

 In particular, letting $\omega = 0$, we see that $H = (H \cap A)
\semiprod (H \cap N)$.
 \end{lem}

\begin{lem}[{cf.~\cite[pf.~of Thm.~3.2.5, p.~42]{Zimmer}}]
\label{Zar-norm}
 If $H$ is a closed, connected subgroup of~$G$, then the Zariski closure
of~$H$ normalizes~$H$.
 \end{lem}

\begin{lem}[{cf.~\cite[Thm.~10.6, pp.~137--138]{Borel-AlgicGrp}}]
\label{algic->semiprod}
 If $H$ is a closed connected subgroup of~$AN$, such that $H$
has finite index in its Zariski closure, then $H$ can be written as a
semidirect product $H = T \semiprod U$, where $U$ is a subgroup of~$N$ and
$T$ is conjugate, via an element of~$N$, to a subgroup of~$A$.
 \end{lem}

\begin{lem} \label{HN=AN}
 Let $H$ be a closed, connected subgroup of~$AN$.
 If $\dim H - \dim(H \cap N) \ge \Rrank G$, then $H$ contains a
conjugate of~$A$, so $H$ is a Cartan-decomposition subgroup.
 \end{lem}

\begin{proof}
 Let $\pi  \colon AN \to A$ be the projection with kernel~$N$. Since
 $$ \Rrank G \le (\dim H) - \dim(H \cap N) = \dim\bigl(\pi(H) \bigr) \le
\dim A = \Rrank G ,$$
 we must have $\pi(H) = A$, so $HN = AN$. Therefore, letting
$\closure{H}$ be the Zariski closure of~$H$, we may assume
that~$\closure{H}$ contains~$A$ \see{algic->semiprod}, by replacing~$H$
with a conjugate subgroup. So $H$ is normalized by~$A$ \see{Zar-norm}.
Therefore, since $\pi(H) = A$, we conclude that $A \subset H$
\see{normedbyA}.
 \end{proof}

All maximal compact subgroups of any connected Lie group are conjugate
\cite[Thm.~XV.3.1, p.~180--181]{Hochschild-Lie}, so the quantity $\dim H
- \dim K_H$ in the statement of the following lemma is independent of
the choice of~$K_H$.

\begin{lem} \label{HcanbeAN}
 Let $H$ be a closed, connected subgroup of~$G$. Then there is a closed,
connected subgroup~$H'$ of~$G$ and a compact subgroup~$C$ of~$G$, such
that $CH = CH'$, and $H'$ is conjugate to a subgroup of~$AN$.
Furthermore, there is a continuous function $f \colon H \to H'$ with
$f(h) \in Ch$ for every~$h \in H$, and we have $\dim H' =
 \dim H - \dim K_H$, where $K_H$ is a maximal compact subgroup of~$H$.
 \end{lem}

\begin{proof}[Sketch of proof]
 Let $L$ be a maximal connected semisimple subgroup of~$H$, and let $T$
be a maximal compact torus of the Zariski closure of $\Rad H$, the
solvable radical of~$H$. Replacing~$T$ by a conjugate torus, we may
assume that $L$ centralizes~$T$. Let $L = K_L A_L N_L$ be the Iwasawa
decomposition of~$L$. From L.~Auslander's nilshadow construction
(cf.~\cite[\S4]{Witte}), we know that there is a unique
connected, closed subgroup~$R$ of~$G$, such that $R$ is conjugate to a
subgroup of~$AN$, and $RT = (\Rad H)T$. Let $H' = A_L N_L R$, and note
that the uniqueness of~$R$ implies that $H$ normalizes~$R$, so $H'$ is a
subgroup of~$G$. Then $K_L T H' = TH = K_L T H$.

Define $f \colon H \to H'$ by specifying that $h \in K_LT \cdot f(h)$.
Because $H'$ is conjugate to a subgroup of $AN$, which has no nontrivial
compact subgroups, we have $(K_L T) \cap H' = e$, so $f(h)$ is well
defined.
 \end{proof}

\begin{rem} \label{CDS-onlyAN}
 The proof of Lemma~\ref{HcanbeAN} is constructive. Furthermore, given a
subgroup~$H$ and the corresponding subgroup~$H'$, it clear that $H$ is a
Cartan-decomposition subgroup if and only if $H'$ is a
Cartan-decomposition subgroup. Thus, to characterize all the
Cartan-decomposition subgroups of~$G$, it suffices to find all the
Cartan-decomposition subgroups that are contained in~$AN$.
 \end{rem}

\begin{rem}
 Our restriction to closed subgroups in the definition of
Cartan-decomposition subgroups is not very important. Namely, if one were
to allow non-closed subgroups, then one would prove that a subgroup is a
Cartan-decomposition subgroup if and only if its closure is a
Cartan-decomposition subgroup. This follows from the theorem of M.~Goto
and A.~Malcev (independently) that if $H$ is a connected Lie subgroup
of~$G$, then there is a compact subgroup~$C$ of~$G$, such that $CH$ is
the closure of~$H$ \cite[Thm.~1.3]{Poguntke}. (This theorem can be
derived from the proof of Lemma~\ref{HcanbeAN}, because every connected
subgroup of~$AN$ is closed.)
 \end{rem}

\section{General results on Cartan-decomposition subgroups}
\label{CDS-general}

\begin{notation}
 We employ the usual Big Oh and little oh notation: for functions
$f_1,f_2$ on~$H$, and a subset~$Z$ of~$H$, we say \emph{$f_1 = O(f_2)$
for $z \in Z$} if there is a constant~$C$, such that, for all large $z
\in Z$, we have $\|f_1(z)\| \le C \|f_2(z)\|$. (The values of each~$f_i$
are assumed to belong to some finite-dimensional normed vector space,
typically either~$\real$ or a space of real matrices. Which particular
norm is used does not matter, because all norms are equivalent up to a
bounded factor.) We say \emph{$f_1 = o(f_2)$ for $z \in Z$} if
$\|f_1(z)\|/\|f_2(z)\| \to 0$ as $z \to \infty$. Also, we write $f_1
\asymp f_2$ if $f_1 = O(f_2)$ and $f_2 = O(f_1)$.
 \end{notation}

 \begin{notation} \label{fund-reps}
 It is known (cf.~\cite[Lem.~2.3]{Benoist}) that there exist irreducible
real finite-dimensional representations $\rho _i$, $i=1, \ldots , k$ of
(a finite cover of)~$G$, such that the highest weight space of
each~$\rho_i$ is one-dimensional and if $\chi _i$ is the highest weight
of $\rho _i$, then $\{\, \chi _i \mid i=1, \cdots ,k \,\}$ is a basis of
the vector space $A^*$ of all continuous group homomorphisms from~$A$
to~$\real^+$ (where the vector-space structure on~$A^*$ is defined by $(s
\alpha + t \beta)(a) = \alpha(a)^s \beta(a)^t$ for $s,t \in \real$,
$\alpha,\beta \in A^*$, and $a \in A$). In particular, we have $k = \Rrank
G$.

When we have fixed a particular choice of $\rho_i$, $i=1, \ldots , k$,
we may refer to $\rho_1,\ldots,\rho_k$ as the \emph{fundamental}
representations of~$G$.
 \end{notation}

\begin{prop}[{Benoist \cite[Lem.~2.4]{Benoist}}] \label{norm}
 For each $i = 1,\ldots, \Rrank G$, we have
 $\chi _i \bigl(\mu (g) \bigr) \asymp \rho _i (g)$, for $g\in G$.
\end{prop}

Because $\mu(g)$ is determined by the values
$\chi_1\bigl(\mu(g)\bigr),\ldots,\chi_k\bigl(\mu(g)\bigr)$ (and
$\chi_i(a) = |\chi_i(a)|$ for each $a \in A$), it follows from the
preceding proposition that the Cartan projection $\mu (g)$ can be
calculated with bounded error by finding the norms of $\rho _1
(g),\ldots,\rho_k(g)$. (The error bound depends only on~$G$; it is
independent of~$g$). This is theoretically useful (see, for example,
Corollaries~\ref{bdd-change}, \ref{proper}, and~\ref{conj-alm-same}, and
note that Theorem~\ref{CDSvsmu} is a special case of
Corollary~\ref{bdd-change}) and is also the method we use in practice in
Sections~\ref{SL3R}, \ref{SO2n-HinN-sect} and~\ref{SO2n-notinN-sect} to
calculate the image under~$\mu$ of subgroups of $\SL(3,\real)$ and
$\SO(2,n)$.

\begin{cor}[{Benoist \cite[Prop.~5.1]{Benoist}}] \label{bdd-change}
 For any compact set~$C_1$ in~$G$, there is a compact set~$C_2$ in~$A$,
such that $\mu(C_1 g C_1) \subset \mu(g) C_2$, for all $g \in G$. 
 \end{cor}

\begin{cor}[{Benoist \cite[Prop.~1.5]{Benoist}, Kobayashi
\cite[Cor.~3.5]{Kobayashi-criterion}}] \label{proper}
 Let $H_1$ and~$H_2$ be closed subgroups of~$G$. The subgroup~$H_1$ acts
properly on~$H_2$ if and only if, for every compact subset~$C$ of~$A$,
the intersection $\bigl( \mu(H_1) C \bigr) \cap \mu(H_2)$ is compact.
 \end{cor}

\begin{notation} 
 For subsets $X,Y \subset A^+$, we write $X \approx Y$ if there is a
compact subset $C$ of~$A$ with $X \subset YC$ and $Y \subset XC$.
 \end{notation}

It is obvious from the definition that if $H$ is a Cartan-decomposition
subgroup, then every conjugate of~$H$ is also a Cartan-decomposition
subgroup. In other words, if $\mu(H) \approx A^+$, then $\mu(g^{-1} H g)
\approx A^+$. The following corollary is a generalization of this
observation.

\begin{cor} \label{conj-alm-same}
 If $H$ is a subgroup of~$G$, then $\mu(g^{-1} H g) \approx \mu(H)$, for
every $g \in G$.

In particular, every conjugate of a Cartan-decomposition subgroup is a
Cartan-decomposition subgroup.
 \end{cor}

\begin{notation} \label{k1k2}
 Suppose that $\Rrank G=2$ and let $\{\alpha _1 ,\alpha _2\}$ be the set
of simple roots with respect to $A^+$ in $G$. We denote by $L_i$ the wall
of $A^+$ defined by $\alpha _i=1$ for each $i=1,2$. Since, by definition,
$\{\chi _1, \chi _2\}$ is a basis of~$A^*$, there exists some real
number~$k_i$ such that
 $$ L_i=\{\, a\in A^+ \mid \chi _1 (a)^{k_i}= \chi _2 (a) \, \} .$$
 (Note that $\alpha_i$ cannot be a scalar multiple of~$\chi_1$, because
$\chi_1$, being a highest weight, is in~$A^+$, but $\alpha$, being a
simple root, is not.) Although we do not need this general fact, we
mention that $k_i$ is always a rational number.
 \end{notation}

We now introduce convenient notation for describing the image of a
subgroup under the Cartan projection~$\mu$.

\begin{notation} \label{approx}
 Suppose that $\rho$ is a representation of~$G$ and that $G$ is a matrix
group (that is, $G \subset \GL(\ell,\real)$, for some~$\ell$), and assume
that the two fundamental representations of~$G$ are $\rho_1(g) = g$ and
$\rho_2(g) = \rho(g)$ (cf.~Notation~\ref{fund-reps}).
 For functions $f_1,f_2 \colon \real^+ \to \real^+$, and a subgroup~$H$
of~$G$, we write 
 $\mu(H) \approx \muH{f_1(\|h\|)}{f_2(\|h\|)}$ if, for every
sufficiently large $C > 1$, we have $$ \mu(H) \approx \bigset{a \in A^+}
 { C^{-1} f_1 \bigl( \|a\| \bigr) \le \|\rho(a)\| \le C f_2 \bigl( \|a\|
\bigr)} .$$
 (If $f_1$~and~$f_2$ are monomials, or other very tame functions, then
it does not matter which particular norm is used.)
 In particular, from Proposition~\ref{norm}, we see that $H$ is a
Cartan-decomposition subgroup of~$G$ if and only if $\mu(H) \approx
\muH{\|h\|^{k_1}}{\|h\|^{k_2}}$, where $k_1$~and~$k_2$ are as described
in Notation~\ref{k1k2}.
 \end{notation}

If $\Rrank G = 2$, Proposition~\ref{CDS<>seq} provides a simple way to
determine whether or not a subgroup~$H$ of~$AN$ is a Cartan-decomposition
subgroup. One simply needs to determine whether or not there are
arbitrarily large elements $h_1$ and~$h_2$ of~$H$, such that
$\|\rho_2(h_i)\|$ is approximately $\|\rho_1(h_i)\|^{k_i}$. 

The proof of Proposition~\ref{CDS<>seq} uses some basic properties of
polynomials of exponentials. The class of $\exp$-definable functions
\see{definable} is much more general, and the theory of the real
exponential function is o-minimal \cite{WilkiesJAMS, WilkiesICM},
so general properties of o-minimal structures
\cite{PillaySteinhorn, KnightPillaySteinhorn, vandenDries} tell us that
every $\exp$-definable function behaves very much like an ordinary
polynomial.

The following definition assumes some familiarity with first-order logic.

\begin{defn}[{cf.~\cite[(5.3)]{vandenDries}}] \label{definable}
 A subset~$X$ of~$\real^n$ is \emph{$\exp$-definable} if there are  real
numbers $a_1,\ldots,a_m$, and a first-order sentence
$\phi(x_1,\ldots,x_n,y_1,\ldots,y_m)$ in the language of ordered fields
augmented with the exponential function~$\exp$, such that
 $$ X = \bigl\{\, (x_1,\ldots,x_n) \in \real^n \mid
 \phi(x_1,\ldots,x_n,a_1,\ldots,a_m) \,\bigr\} .$$
 A map $f \colon X \to Y$ is \emph{$\exp$-definable} if its graph is an
$\exp$-definable subset of $X \times Y$.
 \end{defn}

We remark that every $\exp$-definable set~$X$ has a well-defined
dimension \cite[p.~5 and p.~63]{vandenDries}. Namely, $X$ can be written
as the union of finitely many cells, and the dimension of~$X$ is the
maximum of the dimensions of these cells.

\begin{lem} \label{mu(H)definable}
 If $H$ is a closed, connected subgroup of~$AN$, then $\mu(H)$ is an
$\exp$-definable subset of~$A$, and $\dim \mu(H) \le \dim H$.
 \end{lem}

\begin{proof}
 Because $H = \exp \Lie h$, it is easy to see that $H$ is an
$\exp$-definable subset of~$G$. Furthermore, the map~$\mu$ is obviously
$\exp$-definable. Therefore, $\mu(H)$ is an $\exp$-definable subset
of~$A$, and $\dim \mu(H) \le \dim H$ \cite[Cor.~4.1.6(ii),
p.~66]{vandenDries}.
 \end{proof}

\begin{lem}[{\cite[Props.~3.2.8 and~6.3.2, pp.~57, 100]{vandenDries}}]
\label{def->fincomps}
 Any $\exp$-definable set has only finitely many path-connected
components.
 \end{lem}

\begin{prop} \label{dim>Rrank}
 Let $H$ be a closed, connected subgroup of~$G$, and let $K_H$ be a
maximal compact subgroup of~$H$.
 If $H$ is a Cartan-decomposition subgroup of~$G$, or, more generally, if
no closed, noncompact subgroup of~$G$ acts properly on $G/H$, then $\dim
H - \dim K_H \ge \Rrank G$.
 \end{prop}

\begin{proof}
 We may assume that $H \subset AN$ \see{HcanbeAN}, so $\dim H - \dim K_H
= \dim H$. We prove the contrapositive: suppose that $\dim H < \Rrank G$.
Let $\Lie s = \{\, a \in \Lie a^+ \mid \exp(a) \in \mu(H) \,\}$. Then
$\Lie s$ is $\exp$-definable and $\dim \Lie s \le \dim H$ \see{mu(H)definable}.
Let $\Lie a_\infty$ be the sphere at infinity, and let $\Lie s_\infty$ be
the set of accumulation points of~$\Lie s$ in~$\Lie a_\infty$. For any
nonempty $\exp$-definable set~$S$, we have $\dim \bdry S < \dim S$, where $\bdry
S$ is the complement of~$S$ in its closure \cite[Thm.~4.1.8,
p.~67]{vandenDries}, so we know that $\dim \Lie s_\infty < \dim \Lie s$.
Therefore, $\dim \Lie s_\infty < \dim \Lie a - 1 = \dim \Lie a_\infty$,
so $\Lie s_\infty$ does not contain any nonempty open subset of~$\Lie
a_\infty$. On the other hand, it is obvious that $\Lie a^+_\infty$ does
contain an open subset of~$\Lie a_\infty$, so we conclude that some point
$r \in \Lie a_\infty$ does not belong to~$\Lie s_\infty$. Let $\Lie r^+$
be the ray in~$\Lie a^+$ from~$0$ in the direction~$r$. Because $r \notin
\Lie s_\infty$, we know that, for every compact subset~$\Lie c$ of~$\Lie
a$, the intersection $\Lie r^+ \cap (\Lie s + \Lie c)$ is bounded. Thus,
Corollary~\ref{bdd-change} implies that, for every compact subset~$C$
of~$G$, the intersection $\exp(\Lie r^+) \cap CHC$ is bounded. Therefore,
the one-parameter subgroup $\exp(\Lie r^+ \cup -\Lie r^+)$ acts properly
on~$G/H$, so $H$ is not a Cartan-decomposition subgroup
\see{Calabi-Markus}. 
 \end{proof}

The converse of Proposition~\ref{Calabi-Markus} is not known to hold in
general. However, the following proposition (combined with
Proposition~\ref{Rrank1-CDS}) shows that the converse does hold if
$\Rrank G \le 2$. Condition~\fullref{CDS<>seq}{-seq} is used throughout
Sections~\ref{SL3R}, \ref{SO2n-HinN-sect}, and~\ref{SO2n-notinN-sect} to
determine whether or not a subgroup of~$G$ is a Cartan-decomposition
subgroup.

\begin{prop} \label{CDS<>seq}
 Assume that $\Rrank G = 2$. If $H$ is a connected subgroup of~$AN$, then
the following are equivalent:
 \begin{enumerate}
 \item \label{CDS<>seq-CDS}
 $H$ is a Cartan-decomposition subgroup of~$G$.
 \item \label{CDS<>seq-CM}
 No closed, noncompact subgroup of~$G$ acts properly on $G/H$.
 \item \label{CDS<>seq-seq}
 We have $\dim H \ge 2$ and, for each $i=1,2$, there exists a sequence
$h_i(n)$ in $H$, such that $h_i (n) \to \infty$ as $n \to \infty$, and
 $\rho _2 \bigl(h_i (n)\bigr) \asymp \|\rho _1
\bigl(h_i(n)\bigr)\|^{k_i}$ for $n \in \integer^+$.
 \end{enumerate}
 \end{prop}

\begin{proof}
 (\ref{CDS<>seq-CDS} $\Rightarrow$ \ref{CDS<>seq-CM}) See
Proposition~\ref{Calabi-Markus}.

 (\ref{CDS<>seq-CM} $\Rightarrow$ \ref{CDS<>seq-seq}) This follows from
Proposition~\ref{dim>Rrank}, Theorem~\ref{CDSvsmu}, and
Proposition~\ref{norm}.

 (\ref{CDS<>seq-seq} $\Rightarrow$ \ref{CDS<>seq-CDS})
 We begin by showing that, instead of only sequences $h_1(n)$
and~$h_2(n)$, there are two continuous curves $h_1'(t)$ and $h_2'(t)$, $t
\in [0, \infty)$, in~$H$, such that $h_i' (t) \to \infty$ as $t \to
\infty$, and, for each $i = 1,2$, we have
 $\rho _2 \bigl(h_i' (t) \bigr) \asymp \|\rho _1
\bigl(h_i'(t)\bigr)\|^{k_i}$, for $t \in [0, \infty)$. Fix $i \in
\{1,2\}$. For each~$n$, and each $C > 0$ let 
 $$A(n,C) = \{\, a \in \mu(H) \mid \|a\| > n \mbox{ and }
 C^{-1} \|\rho _1 (a)\|^{k_i}
 < \|\rho _2 (a) \|
 < C \|\rho _1 (a)\|^{k_i} \,\} .$$
 The existence of the sequence~$\{h_i(n)\}$ implies that there is some
$C_0 > 0$, such that, for every~$n$, the set $A(n,C_0)$ is nonempty.
 Because $A(n,C_0)$ is $\exp$-definable \cf{mu(H)definable},
Lemma~\ref{def->fincomps} implies that we may choose an unbounded
path-connected component~$A_n$ of $A(n,C_0)$, for each~$n$. Furthermore,
we may assume that $A_{n+1} \subset A_n$. Now, for each natural
number~$n$, choose a point~$a_n \in A_n$, and let $h_i'(t)_{n \le t \le
n+1}$ be any path in~$A_n$ from~$a_n$ to~$a_{n+1}$.
 
Since $H$ is homeomorphic to some Euclidean space~$\real ^m$, with $m \ge
2$, it is easy to find a continuous and proper map $\Phi\colon [1,2]
\times \real ^+ \to H$ such that $\Phi (i,t)= h_i'(t)$ for $i=1,2$ and
for all $t \in \real^+$. From Proposition~\ref{norm} and the definition
of~$k_i$, we know that the curve $\mu\bigl(h_i'(t) \bigr)$ stays within a
bounded distance from the wall~$L_i$; say $\operatorname{dist}\bigl[
\bigl(\Phi (i,t) \bigr), L_i \bigr] < C$ for all~$t$. We may assume that
$C$ is large enough that $\operatorname{dist}\bigl( \Phi(s,1), e \bigr) <
C$ for all $s \in [1,2]$. Then an elementary homotopy argument shows that
$\mu \bigl[ \Phi\bigl( [1,2] \times \real ^+ \bigr) \bigr]$ contains
 $\{\, a \in A^+ \mid \operatorname{dist}(a, L_1 \cup L_2) > C \, \}$, 
 so $\mu \bigl[ \Phi\bigl( [1,2] \times \real ^+ \bigr) \bigr] \approx
A^+$. Because $\mu(H) \supset \mu \bigl[ \Phi\bigl( [1,2] \times \real ^+
\bigr) \bigr]$, we conclude from Theorem~\ref{CDSvsmu} that $H$ is a
Cartan-decomposition subgroup.
 \end{proof}

\begin{rem}
 When verifying Condition~\fullref{CDS<>seq}{-seq} for a specific example,
one may use any matrix norm, because any two norms are equivalent up to a
bounded factor. In practice, the authors use the maximum absolute value
of the matrix entries, but the reader is free to make another choice, 
because we calculate our results only to within a bounded factor.  In our
applications to $\SL(3,\real)$ and $\SO(2,n)$, we have $\rho_2 = \rho_1
\wedge \rho_1$, the second exterior power of~$\rho_1$. So we calculate
$\|\rho_2(h)\|$ as the maximum  absolute value among the determinants of
all the $2 \times 2$ submatrices of~$\rho_1(h)$.
 \end{rem}

If $\dim H = 1$ and $\Rrank G \ge 2$, then we
know, from Proposition~\ref{dim>Rrank}, that $H$ is not a
Cartan-decomposition subgroup of~$G$. For completeness, the following
simple proposition describes $\mu(H)$ fairly explicitly. Namely,
$\mu(H)$ consists of either one or two rays, or one or two logarithmic
curves.

\begin{prop} \label{dim1-mu(H)}
 Assume that $\Rrank(G) \ge 2$, and let $H$ be a nontrivial one-parameter
subgroup of~$AN$.
 \begin{enumerate}
 \item \label{dim1-mu(H)-A}
 If $H$ is conjugate to a subgroup of~$A$, then there is a
ray~$R$ in~$A^+$, such that $\mu(H) \approx R \cup i(R)$, where $i \colon
A^+ \to A^+$ is the opposition involution \see{oppinv}.
 \item \label{dim1-mu(H)-N}
 If $H \subset N$, then there is a
ray~$R$ in~$A^+$, such that $\mu(H) \approx R$.
 \item \label{dim1-mu(H)-neither}
 If neither \pref{dim1-mu(H)-A} nor~\pref{dim1-mu(H)-N} applies,
then there is a ray~$R$ in~$A^+$, a ray~$R'$ in~$A$ that is
perpendicular to~$R$, and a positive number~$k$, such that  $\mu(H)
\approx X \cup i(X)$, where 
 $X = \{\, rs \mid r \in R, ~ s \in R', ~ \|s\| = (\log\|r\|)^k \,\}$.
 \end{enumerate}
 \end{prop}

\begin{proof}
 Let $H = \{h^t\}$. Replacing $H$ by a conjugate, we may assume that $h^t
= a^t u^t$, where $\{a^t\}$ is a one-parameter subgroup of~$A$ and
$\{u^t\}$ is a one-parameter subgroup of~$N$, such that $\{a^t\}$
centralizes~$\{u^t\}$ \cf{not-semi}. The
subgroup $\{u^t\}$ is contained in a subgroup~$S$ of~$G$ that is locally
isomorphic to~$\SL(2,\real)$ \cite[Thm.~17(1), p.~100]{Jacobson}.
Replacing $H$ by a conjugate, we may assume that the intersection $A_S =
A \cap S$ is nontrivial.

 \pref{dim1-mu(H)-A} We have $\{h^t\} = \{a^t\}$. Replacing $H$ by a
conjugate under the Weyl group, we may assume that $a^t \in A^+$ for all
$t \ge 0$. Let $R = \{\, a^t \mid t \ge 0 \,\}$.

 \pref{dim1-mu(H)-N} We have $\{h^t\} = \{u^t\}$, so $H$ is contained
in~$S$. Then, because both $H$ and~$A_S$ are Cartan-decomposition
subgroups of~$S$ \see{Rrank1-CDS}, we have $\mu(H) \asymp \mu(S) \asymp
\mu(A_S^+)$. Let $R = \mu(A_S^+)$.

 \pref{dim1-mu(H)-neither} Both $a^t$ and~$u^t$ are nontrivial. Because
$\{a^t\}$ centralizes~$\{u^t\}$, and hence centralizes all of~$S$, we
know that $A_S$ is perpendicular to~$\{a^t\}$ (with respect to the
Killing form). Let $\mu_S \colon S \to A_S^+$ be the Cartan projection
of~$S$. (We may assume that $K \cap S$ is a maximal compact subgroup
of~$S$.) Assuming, for simplicity, that $\{a^t \mid t \ge 0\} \subset
A^+$, we see from the proof of Proposition~\ref{rootsemi-CDS} that
$\mu(h^t) = a^t \mu_S(u^t)$. Then, because $\|a^t\|$ grows exponentially
and $\|u^t\|$ grows polynomially, the desired conclusion follows, with $R
= \{a^t \mid t \ge 0\}$ and $R' = A_S^+$.
 \end{proof}

\begin{prop} \label{rootsemi-CDS}
 Assume that $\mathop{\real\text{\rm-rank}}(G) = 2$, and that
$\omega$ is a real root of~$G$. Let $H$ be a closed, connected subgroup
of~$G$, such that $\Lie h = \Lie t + \Lie u$, for some one-dimensional
subspace~$\Lie t$ of~$\Lie a$ and some nontrivial subalgebra~$\Lie u$ of
$\Lie u_\omega + \Lie u_{2\omega}$.
 \begin{enumerate}
 \item Let $\Lie a_\omega = [\Lie u_\omega, \Lie u_{-\omega}] \cap \Lie
a$, a one-dimensional subspace of~$\Lie a$.
 \item Choose rays~$\Lie t^+$ in~$\Lie t$ and~$\Lie a_\omega^+$ in~$\Lie
a_\omega$, such that $\langle t \mid a \rangle \ge 0$ for every $t \in
\Lie t^+$ and every $a \in \Lie a_\omega^+$, where $\langle \cdot \mid
\cdot \rangle$ is the inner product on~$\Lie a$ defined by the Killing
form.
 \end{enumerate}
 The subgroup~$H$ is a Cartan-decomposition subgroup if and only if $\Lie
t^+ - \{0\}$ is {\bf not} contained  in the interior of the
region~$\chamber$ defined as follows:
 \begin{enumerate} \renewcommand{\theenumi}{\alph{enumi}}
 \item if $\Lie a_\omega^+$ is contained in the interior of a Weyl
chamber, let $\chamber$ be the interior of that Weyl chamber; 
 \item if $\Lie a_\omega^+$ is the wall between two Weyl chambers
$\chamber_1$ and~$\chamber_2$, let $\chamber$ be the interior of
$\chamber_1 \cup \Lie a_\omega^+ \cup \chamber_2$.
 \end{enumerate}
 \end{prop}

\begin{proof}
 Let $M = \langle U_\omega ,U_{2\omega }, U_{-\omega} ,U_{-2\omega }
\rangle$, so $M$ is semisimple of real rank one. (The only real roots
of~$M$ are~$\pm \omega$ and, possibly,~$\pm 2\omega$.) Let $A_\omega = M
\cap A$, so $A_\omega$ is a maximal split torus of~$M$, and note that
the Lie algebra of~$A_\omega$ is~$\Lie a_\omega$. The choice of the
ray~$\Lie a_\omega^+$ determines a corresponding Weyl
chamber~$A_\omega^+$ in~$A_\omega$.

We may assume that $K_M = K \cap M$ is a maximal compact subgroup of~$M$,
and let $\mu_M \colon M \to A_\omega^+$ be the corresponding Cartan
projection.

 Let $S = C_A(M)$. The subgroup $MA$ is reductive, so we have $MA = MS$.
In particular, any element of~$A$ may be written uniquely in the form $as$
with $a \in A_\omega$ and $s \in S$. We extend $\mu_M$ to a map
$\mu_{MA} \colon MA \to A$ by defining $\mu_{MA}(ms) = \mu_M(m) s$ for $m
\in M$ and $s \in S$.

Because $H \cap N \subset U_\omega U_{2\omega} \subset M$ and $T \subset
A$, we have $H \subset MA$.

\begin{claim}
 We have $\mu_{MA}(H) = A_\omega^+ \, \mu_{MA}(T)$.
 \end{claim}
 Given $h \in H$, we may write $h = ut$, with $u \in H \cap N$ and $t \in
T$. Furthermore, we may write $t = as$, with $a \in A_\omega$ and $s \in
S$. Choose $k,k' \in K_M$ with $k ua k' = \mu_M(ua) \in A_\omega^+$. 

Assume for simplicity that $G$ is a matrix group, that is, $G \subset
\GL(\ell,\real)$, for some~$\ell$. Then, for $g \in G$, we may let
$\|g\|^2$ be the sum of the squares of the matrix entries of~$g$. Assume,
furthermore, that $K = G \cap \SO(\ell)$, that $A$ is the group of
diagonal matrices in~$G$, and that $N$ is the group of unipotent
upper-triangular matrices in~$G$. Then it is clear that $\|ua\| \ge
\|a\|$ and, because $\|\cdot\|$ is bi-$K$-invariant, that $\| \mu(g)\| =
\|g\|$ for all $g \in G$.

 Because
 $$ \|\mu_M(ua)\| = \|ua\| \ge \|a\| = \|\mu_M(a)\| ,$$
 and both of $\mu_M(ua)$ and~$\mu_M(a)$ are in~$A_\omega^+$, there is
some $a^+ \in A_\omega^+$ with $\mu_M(ua) = a^+ \mu_M(a)$. Therefore
 $$\mu_{MA}(h) = \mu_{MA}(uas) = \mu_M(ua)s = a^+ \mu_M(a) s = a^+
\mu_{MA}(as) \in A_\omega^+ \mu_{MA}(T) .$$
 Therefore $\mu_{MA}(H) \subset A_\omega^+ \, \mu_{MA}(T)$.

Conversely, given $a^+ \in A_\omega^+$ and $t \in T$, write $t = as$ with
$a \in A_\omega$ and $s \in S$. Because
 $\|a^+ \mu_M(a)\| \ge \|\mu_M(a)\| = \|a\|$,
 there is some $u \in H \cap N$ with $\|ua\| = \|a^+ \mu_M(a)\|$. Then
$\mu_M(ua) = a^+ \mu_M(a)$, because there is only one element
of~$A_\omega^+$ with any given norm, so 
 $$\mu_{MA}(ut) = \mu_M(ua)s = a^+ \mu_M(a) s = a^+
\mu_{MA}(t) .$$
 This completes the proof of the claim.

\medskip

Note that $\mu_{MA}(g) \in K_M g K_M \subset KgK$, so $\mu(g) = \mu
\bigl( \mu_{MA}(g) \bigr)$, for all $g \in MA$.

Now, to clarify the situation, let us define $\mu_{\Lie m + \Lie
a} \colon \Lie m + \Lie a \to \Lie a$ by $\mu_{MA}(\exp z) =
\exp \bigl(\mu_{\Lie m + \Lie a}(z) \bigr)$ for $z \in \Lie m + \Lie
a$, and let us introduce a convenient coordinate system on the Lie
algebra~$\Lie a \cong \real^2$. Let the $x$-axis be~$\Lie s$, the Lie
algebra of~$S$, and let the positive $y$-axis be~$\Lie a_\omega^+$. In
these coordinates, the restriction of~$\mu_{\Lie m + \Lie a}$ to~$\Lie a$
is given by $\mu_{\Lie m + \Lie a}(x,y) = (x,|y|)$. The line~$\Lie t$ has
an equation of the form $y = m_{\Lie t} x$, so $\mu_{\Lie m + \Lie a}(T)$
has the equation $y = |m_{\Lie t} x|$. Thus, the claim asserts that
 $$\mu_{\Lie m + \Lie a}(H) = \{\, (x,y) \mid y \ge |m_{\Lie t} x| \,\}
.$$
 For some constant~$m_{\chamber}$, we have
 $$\chamber = \{\, (x,y) \mid y > |m_{\chamber} x| \,\} .$$
 Thus, $\Lie t^+ - \{0\}$ is in the interior of~$\chamber$ if and only if
$|m_{\Lie t}| > |m_{\chamber}|$.

($\Leftarrow$) If $\Lie t^+ - \{0\}$ is not in the interior of~$\chamber$,
then, because $|m_{\Lie t}| \le |m_{\chamber}|$, we see that $\mu_{\Lie m
+ \Lie a}(H)$ contains the closure of~$\chamber$. By definition, we know
that the closure of~$\chamber$ contains a Weyl chamber, so we conclude
that $\mu(H) \supset A^+$, so $H$ is a Cartan-decomposition subgroup.

($\Rightarrow$)  If $\Lie t^+ - \{0\}$ is in the interior of~$\chamber$,
then, because $|m_{\Lie t}| < |m_{\chamber}|$, we see that $\mu_{\Lie m +
\Lie a}(H)$ is a proper subcone of~$\chamber$. If $\chamber$ is a single
Weyl chamber, then this immediately implies that $H$ is not a
Cartan-decomposition subgroup. Now assume the other possibility, namely,
that the $y$-axis~$\Lie a_\omega$ is the wall between two Weyl chambers
$\chamber_1$ and~$\chamber_2$. Assume for simplicity that $\chamber_2 =
A^+$. The reflection across the $y$-axis is the Weyl reflection that maps
$\chamber_1$ to~$\chamber_2$. Because $\mu_{\Lie m + \Lie a}(H)$ is
invariant under this reflection (that is, the inequality defining
$\mu_{\Lie m + \Lie a}(H)$ depends only on~$|x|$, not on~$x$ itself), we
see immediately that $\mu(H) = \chamber_2 \cap \mu_{\Lie m +
\Lie a}(H)$ is a proper subcone of~$\chamber_2$. So $H$ is not a
Cartan-decomposition subgroup.
 \end{proof}

\begin{cor}[of proof] \label{rootsemi-proj}
 Let $G$, $\omega$, $H$, $T$, $M$, and $\Lie a_\omega^+$ be as in the
statement and proof of Proposition~\ref{rootsemi-CDS}. 
 If $H$ is not a Cartan-decomposition subgroup, then the
Cartan projection $\mu(H)$ is the image under a Weyl element of either
the closed, convex cone bounded by~$\mu_{MA}(T)$ {\upshape(}if $\Lie
a_\omega^+$ is in the interior of a Weyl chamber{\upshape)} or the
closed, convex cone bounded by~$A_\omega^+$ and a ray of~$T$
{\upshape(}if $\Lie a_\omega^+$ is a wall of a Weyl chamber{\upshape)}.
 \end{cor}

\begin{eg} \label{SL3-rootsemi-eg}
 Suppose that $H$ is of the form \pref{SL3-rootsemi} for some $p,q \in
\real$, not both zero. As an illustration of the application of
Proposition~\ref{rootsemi-CDS}, we show that $H$ is a
Cartan-decomposition subgroup of $\SL(3,\real)$ if and only if either $p+q
\le -\max\{p,q\}$ or $p+q \ge - \min\{p,q\}$. Let
 $\Lie a_\omega =  \{\, \diag(a,-a,0) \mid a \ge 0 \,\}$.
 Replacing $p$~and~$q$ by their negatives results in the same group, so
we may assume that $p \ge q$. Then we may choose
 $\Lie t^+ = \{\, \diag(pt,qt,-(p+q)t \mid t \ge 0 \,\}$.
 The Weyl chamber 
 $$\chamber = \bigl\{ \, \diag \bigl( a,b,-(a+b) \bigr) \mid a \ge -(a+b)
\ge b \, \bigr\}$$
 contains $\Lie a_\omega^+$ in its interior. It is clear that $\Lie
t^+$ is in the interior of~$\chamber$ if and only if $p > -(p+q) > q$. So
$\Lie t^+$ is a Cartan-decomposition subgroup if and only if this
condition fails, that is, if and only if either $-(p+q) \ge p$ or $-(p+q)
\le q$. Because we have assumed that $p \ge q$, which means that $p =
\max\{p,q\}$ and $q = \min\{p,q\}$, this is the desired conclusion. (The
condition in Theorem~\ref{SL3-minl} that $\max\{p,q\} = 1$ and
$\min\{p,q\} \ge -1/2$ is obtained by assuming that $\max\{p,q\} = 1$,
instead of assuming that $p \ge q$, as we have here. Replacing the pair
$(p,q)$ with any nonzero scalar multiple results in the same subgroup,
so assuming that $\max\{p,q\} = 1$ results in no loss of generality.)
 \end{eg}

 The following proposition shows that characterizing the
Cartan-decomposition subgroups of a reductive group reduces to the
problem of characterizing the Cartan-decomposition subgroups of its
almost simple factors. Thus, our standing assumption that $G$ is almost
simple is not as restrictive as it might seem.

\begin{prop} \label{reductive}
 Suppose that $\tilde G$ is a connected, reductive, linear Lie group, let
$\tilde G = \tilde K \tilde A \tilde N$ be an Iwasawa decomposition
of~$\tilde G$, and write $\tilde G = G_1 G_2 \cdots G_m Z$, where $Z$~is
the center of~$\tilde G$, and each $G_i$~is a connected, almost simple,
normal subgroup of~$\tilde G$. Let $H$ be a closed, connected
subgroup of $\tilde A \tilde N$.
 \begin{enumerate}
 \item \label{reductive-CDS}
 $H$ is a Cartan-decomposition subgroup of~$\tilde G$ if and only
if {\upshape(a)}~$H$ contains~$Z \cap \tilde A$, and {\upshape(b)}~for
each~$i$, the intersection $H \cap G_i$ is a Cartan-decomposition
subgroup of~$G_i$.
 \item \label{reductive-CM}
 No closed, noncompact subgroup of~$\tilde G$
acts properly on $\tilde G/H$ if and only if {\upshape(a)}~$H$
contains~$Z \cap \tilde A$, and {\upshape(b)}~for each~$i$, no closed,
noncompact subgroup of~$G_i$ acts properly on $G_i/(H \cap G_i)$.
 \end{enumerate}
 \end{prop}

\begin{proof}
 We prove only the nontrivial direction of each conclusion.

 \pref{reductive-CDS}
 Suppose that $H$ is a Cartan-decomposition subgroup of~$\tilde G$.
(a)~The proof of the claim in Proposition~\ref{rootsemi-CDS}, with~$G_1
G_2 \cdots G_m$ in the role of~$M$, shows that $H$ contains~$Z \cap
\tilde A$. (b)~From Theorem~\ref{CDSvsmu}, we know that there is a
compact set $C \subset A$, such that, for each $a \in A_1^+$, there is
some $h(a) \in H$, such that $\mu \bigl( h(a) \bigr) \in a C$. The Cartan
projection of~$\tilde G/G_1$ is a proper map, and, because $HG_1$ is
closed (see \cite[Thm.~XII.2.2, p.~137]{Hochschild-Lie}), the natural
map $H/(H \cap G_1) \to \tilde G/G_1$ is a proper map. Hence,
there is a compact subset~$C_1$ of~$H$, such that $h(A_1^+) \subset C_1
(H \cap G_1)$. Then Corollary~\ref{bdd-change} implies that there is a
compact subset~$C_2$ of~$A_1$, such that $\mu(H \cap G_1) C_2 \supset
\mu \bigl( h(A_1^+) \bigr) \approx A_1^+$. Therefore, $H \cap G_1$ is a
Cartan-decomposition subgroup of~$G_1$.

\pref{reductive-CM} Proof by contradiction. Assume that (a) and~(b)
hold, and let $L$ be a closed, noncompact subgroup of~$\tilde G$ that
acts properly on $\tilde G/H$. By passing to a subgroup, we may assume
that $L$ is cyclic. Then there is a one-parameter subgroup~$L'$
of~$\tilde G$, such that $L' \approx L$, and, from Lemma~\ref{HcanbeAN},
we may assume that $L' \subset AN$. From~(a), we see that we may assume
that $L' \subset G_1 \cdots G_m$, so we may write $L' = \{l^t\}$ and
$l^t = l_1^t \cdots l_m^t$, where $\{l_i^t\}$ is a one-parameter subgroup
of~$G_i$. Every nontrivial, connected subgroup of $\tilde A \tilde N$ is
closed and noncompact, so, for each~$i$, (b)~implies that either
$\{l_i^t\}$ is trivial or $\{l_i^t\}$ does not act properly on $G_i/(H
\cap G_i)$. Then, from Corollary~\ref{proper}, we see that, for
each~$i$, there is a compact subset~$C_i$ of~$A \cap G_i$, such that 
 $ \{\, t \ge 0 \mid l_i^t \in \mu_i(H \cap G_i) C_i \,\} $
 is unbounded. (By taking $C_i$ to be a large ball, we may assume that
$C_i$ is $\exp$-definable.) Because this set is $\exp$-definable
\cf{mu(H)definable}, Lemma~\ref{def->fincomps} implies that there is some
$T > 0$, such that, for each~$i$ and for all $t \ge T$, we have $l_i^t
\in \mu_i(H \cap G_i) C_i$. Therefore, for all $t \ge T$, we have
 $l^t \in \mu(H) C_1 C_2 \cdots C_m$, so $L' = \{l^t\}$ does not act
properly on $\tilde G/H$. This is a contradiction.
 \end{proof}

\begin{lem} \label{Gconn}
 Let $H$ be a closed, connected subgroup of~$G$. The subgroup~$H$ is a
Cartan-decomposition subgroup of~$G$ if and only if $H$ is a
Cartan-decomposition subgroup of the identity component of~$G$.
 \end{lem}

\begin{proof}
 Let $G^\circ$ be the identity component of~$G$. Because every element
of the Weyl group of~$G$ has a representative in~$G^\circ$
\cite[Cor.~14.6]{BorelTits}, we see that $G$ and~$G^\circ$ have the same
positive Weyl chamber~$A^+$, and the Cartan projection $G^\circ \to A^+$
is the restriction of the Cartan projection $G \to A^+$. Thus, the
desired conclusion is immediate from Theorem~\ref{CDSvsmu}.
 \end{proof}

\section{Cartan-decomposition subgroups of $\SL(3,\real)$}
 \label{SL3R}

In this section, we find all the Cartan-decomposition subgroups of
$\SL(3,\real)$ \see{SL3-CDS}. With Remark~\ref{CDS-onlyAN} in mind, we
restrict our attention to finding the Cartan-decomposition subgroups
contained in~$AN$. 

\begin{notation}
 To provide a convenient way to refer to specific elements of~$AN$, we
define $h \colon (\real^+)^3 \times \real^3 \to \GL(3,\real)$ by
 $$h(a_1,a_2,a_3,u_1,u_2,u_3) =
 \begin{pmatrix}
 a_1 & u_1 & u_3 \\
 0 & a_2 & u_2 \\
 0 & 0 & a_3 \\
 \end{pmatrix}
 .$$
  \end{notation}

\begin{notation}
Let
 \begin{equation} \label{SL3-A}
 A = \{\, a= h(a_1,a_2,a_3,0,0,0) \mid a_i \in \real^+, ~ a_1
a_2 a_3 = 1 \, \},
 \end{equation}
 $$A^+ = \{\, h(a_1,a_2,a_3,0,0,0) \mid a_1 \ge a_2 \ge a_3 > 0, ~ a_1 a_2
a_3 = 1 \, \},$$
 $$N = \{\, h(1,1,1,u,v,w) \mid u,v,w \in \real \, \},$$
 and $K = \SO(3)$.
  \end{notation}

\begin{notation}
 We let $\alpha$ and~$\beta$ be the simple roots of $\SL(3,\real)$,
defined by $\alpha(a) = a_1/a_2$ and $\beta(a) = a_2/a_3$, for an
element~$a$ of~$A$ of the form~\pref{SL3-A}. Thus,
 \begin{itemize}
 \item the root space $\Lie u_\alpha$ consists of the matrices in which all
entries except~$h_{1,2}$ are~$0$;
 \item the root space $\Lie u_\beta$ consists of the matrices in which all
entries except~$h_{2,3}$ are~$0$; and
 \item the root space $\Lie u_{\alpha+\beta}$ consists of the matrices in
which all entries except~$h_{1,3}$ are~$0$.
 \end{itemize}
 \end{notation}

\begin{thm} \label{SL3-CDS}
 Assume that $G = \SL(3,\real)$. Let $H$ be a closed, connected subgroup
of~$AN$, and assume that $H$ is compatible with~$A$ \see{compatible}.

The
subgroup~$H$ is a Cartan-decomposition subgroup of~$G$ if and only if
either:
 \begin{enumerate}
 \item $\dim H \ge 3$; or
 \item \label{SL3-CDS-A}
 $H = A$; or
 \item \label{SL3-CDS-inN}
 for some $p \neq 0$, we have
 $ H =  \{\, h(1,1,1,s,ps,t) \mid s,t \in \real \, \} \subset N$;
 or
 \item \label{SL3-CDS-alpha}
 for some $p \neq 0$, we have
 $ H =  \{\, h(e^t,e^t,e^{-2t},0,s,ps) \mid s,t \in \real \, \} $;
 or
 \item \label{SL3-CDS-beta}
 for some $p \neq 0$, we have
 $ H =  \{\, h(e^{-2t},e^t,e^t,s,0,ps) \mid s,t \in \real \, \} $;
 or
 \item \label{SL3-CDS-root}
 $H = T \semiprod U_\omega$, for some one-parameter
subgroup~$T$ of~$A$ and some positive root~$\omega$, and $H$ satisfies the
conditions of Proposition~\ref{rootsemi-CDS}.
 \end{enumerate}
 \end{thm}

\begin{rem}
 If $p = 0$ in \fullref{SL3-CDS}{-alpha}
or~\fullref{SL3-CDS}{-beta}, then $H$ is a
Cartan-decomposition subgroup, but $H$ is of the type considered
in~\fullref{SL3-CDS}{-root}. However, if $p=0$ in
\fullref{SL3-CDS}{-inN}, then $H$ is \emph{not} a
Cartan-decomposition subgroup.
 \end{rem}

Theorem~\ref{SL3-CDS} describes the Cartan-decomposition subgroups
of~$\SL(3,\real)$. We now describe the subgroups that are \emph{not}
Cartan-decomposition subgroups.

\begin{cor} \label{SL3notCDS}
 Assume that $G = \SL(3,\real)$. Let $H$ be a closed, connected subgroup
of~$AN$, and assume that $H$ is compatible with~$A$ \see{compatible}.

 The subgroup~$H$ fails to be a Cartan-decomposition subgroup of~$G$ if and
only if either:
 \begin{enumerate}
 \item 
 $\dim H \le 1$; or
 \item \label{SL3notCDS-inN}
 $\Lie h = \Lie u_\alpha + \Lie u_{\alpha+\beta}$ or $\Lie h = \Lie
u_\beta + \Lie u_{\alpha+\beta}$; or
 \item \label{SL3notCDS-semi}
 $\Lie h$ is of the form $\Lie h = \bigl(\ker(\alpha - \beta) \bigr)
+ \Lie u$, for some one-dimensional subspace~$\Lie u$ of $\Lie u_\alpha +
\Lie u_\beta$, with $\Lie u \notin \{\Lie u_\alpha,
\Lie u_\beta \}$; or
 \item \label{SL3notCDS-notsemi}
 $\dim H = 2$ and $H \neq (H \cap A) \semiprod (H \cap N)$; or
 \item $H = T \semiprod U_\omega$, for some one-parameter
subgroup~$T$ of~$A$ and some positive root~$\omega$, and
Proposition~\ref{rootsemi-CDS} implies that $H$ is not a
Cartan-decomposition subgroup.
 \end{enumerate}
\end{cor}

In the course of the proof of Theorem~\ref{SL3-CDS}, we calculate the
Cartan projection of  each subgroup that is not a Cartan-decomposition
subgroup. Proposition~\ref{SL3-Cartanproj} collects these results (see
also Corollary~\ref{rootsemi-proj} and Proposition~\ref{dim1-mu(H)}). The
statement of this proposition and the proof of Theorem~\ref{SL3-CDS} are
based on Proposition~\ref{CDS<>seq}, so we describe the required
representations $\rho_1$ and~$\rho_2$ of~$G$.

\begin{notation} \label{SL3-rho}
 Define a representation 
 $$\rho \colon \SL(3,\real) \to \SL(\real^3 \wedge \real^3)
 \mbox{ by }
 \rho(g) = g \wedge g .$$
 \end{notation}
 In the notation of Proposition~\ref{norm}, we have
 \begin{eqnarray*}
 \rho _1 &=& \mbox {the standard representation on~$\real^3$
 \qquad
($\rho_1(g) = g$)}, \\
 \rho _2 &=& \rho \qquad (= \rho_1 \wedge \rho_1), \\
 k_1 &=& 1/2 \mbox{ \qquad and \qquad } k_2=2 .
 \end{eqnarray*}
 Thus, Proposition~\ref{norm} yields the following
fundamental lemma.

\begin{lem}
 Assume that $G = \SL(3,\real)$.
 A subgroup~$H$ of~$G$ is a Cartan-decomposition subgroup if and only
if $\mu(H) \approx \muH{\|h\|^{1/2}}{\|h\|^2}$, where the
representation~$\rho$ is defined in  Notation~\ref{SL3-rho}.
 \end{lem}

\begin{notation} \label{oppinv}
 Let $i$ be the opposition involution in~$A^+$, that is, for $a\in A^+$,
$i(a)$ is the unique element of~$A^+$ that is conjugate to $a^{-1}$, and
set $B^+=\{\,a\in A^+ \mid i(a)=a \,\}$.
 \end{notation}

\begin{cor} \label{SL3-curves}
 Assume that $G = \SL(3,\real)$.
 Let $H$ be a closed, connected subgroup of~$G$ with $\dim H  - \dim
K_H\ge 2$, where $K_H$ is a maximal compact subgroup of~$H$. The subgroup
 $H$ is a Cartan-decomposition subgroup if and only if there is a
sequence $h_n \to \infty$ in~$H$ with $\rho(h_n) \asymp \|h_n\|^2$.
 \end{cor}

\begin{proof}
 ($\Leftarrow$)
 We may assume that $H \subset AN$ \see{HcanbeAN}. 
If we identify the Lie algebra~$\Lie a$ of~$A$ with the
connected component of~$A$ containing~$e$, then~$A^+$
is a convex cone in~$\Lie a$ and the opposition
involution~$i$ is the reflection in~$A^+$ across the
ray~$B^+$. If $L_1$ and~$L_2$ are the two walls of the
Weyl chamber~$A^+$, then $i(L_1) = L_2$. Therefore, because $\mu(h_n)$ is
a bounded distance from one of the walls, we know that $i\bigl( \mu(h_n)
\bigr)$ is a bounded distance from the other wall. That is,
$\rho\bigl[i\bigl( \mu(h_n) \bigr)\bigr] \asymp \bigl\| i\bigl( \mu(h_n)
\bigr) \bigr\|^{1/2}$. In other words, we have $\rho(h_n^{-1}) \asymp
\|h_n^{-1}\|^{1/2}$. Therefore, using the sequences $h_n$ and~$h_n^{-1}$,
we conclude from Proposition~\ref{CDS<>seq} that $H$ is a
Cartan-decomposition subgroup.
 \end{proof}

We now describe the Cartan projections of the subgroups that are not
Cartan-decomposition subgroups.

\begin{prop} \label{SL3-Cartanproj}
 Assume that $G = \SL(3,\real)$.
 \begin{enumerate}
 \item If $H$ is of type~\fullref{SL3notCDS}{-inN}, then
$\rho(h) \asymp h$ for every $h \in H$.
 \item If $H$ is of type~\fullref{SL3notCDS}{-semi}, then
$\rho(h) \asymp h$ for every $h \in H$.
 \item \label{SL3-Cartanproj-notsemi}
 If $H$ is of type~\fullref{SL3notCDS}{-notsemi}, then
$\mu(H) \asymp \muH{(\|h\| \log \|h\|)^{1/2}}{\|h\|^2/(\log \|h\|)}$.
 \end{enumerate}
 \end{prop}

\begin{proof}[{\bf Proof of Theorem~\ref{SL3-CDS}.}]
 If $A \subset H$, then $H$ is a Cartan-decomposition subgroup (because
$A$ is a Cartan-decomposition subgroup) and we have either $\dim H \ge 3$
or $H = A$. Thus, we may henceforth assume that $A \not\subset H$. Then,
from the proof of Lemma~\ref{HN=AN}, we see that $A$ is not
contained in the Zariski closure~$\closure{H}$ of~$H$. Therefore $\dim(A
\cap \closure{H}) \le 1$, so there are $a,b,c \in \real$ with  
 \begin{equation} \label{SL3-abc}
 \closure{H} \cap A = 
 \{\, h(e^{at}, e^{bt},e^{ct},0,0,0) \mid t \in \real \, \} .
\end{equation}
 and $a+b+c=0$. (If $A \cap \closure{H} = e$, then $a=b=c=0$.) Because a
$1$-dimensional subgroup cannot be a Cartan-decomposition subgroup
\see{dim>Rrank}, we may assume that $\dim H \ge 2$.

\setcounter{case}{0}

\begin{case} \label{HinN}
 Assume that $H \subset N$ and that $\dim H = 2$.
 \end{case}
 We show that $H$ is a Cartan-decomposition subgroup if and only if $H$ is
not normalized by~$A$. (That is, if and only if $H$ is of
type~\pref{SL3-CDS-inN}.) In the case where $H$ is not a
Cartan-decomposition subgroup, we show that $\rho(h) \asymp h$ for every
$h \in H$.

Every $2$-dimensional subalgebra of~$\Lie n$ contains~$\Lie
u_{\alpha+\beta}$, so we must have $\Lie u_{\alpha+\beta} \subset \Lie h$.

($\Rightarrow$) We prove the contrapositive. Thus, we suppose that $H$ is
normalized by~$A$, so $\Lie h$ is a sum of root spaces. Therefore, $\Lie
h$ is either $\Lie u_\alpha + \Lie u_{\alpha+\beta}$ or $\Lie u_{\beta} +
\Lie u_{\alpha+\beta}$. In either case, every element of~$H$ is
conjugate, via~$K$, to an element of~$U_{\alpha+\beta}$, so $\mu(H) =
\mu(U_{\alpha+\beta})$. But $U_{\alpha+\beta}$, being one-dimensional,
cannot be a Cartan-decomposition subgroup \see{dim>Rrank}, so we
conclude that $H$ is not a Cartan-decomposition subgroup.

($\Leftarrow$) Because $H$ is not normalized by~$A$, and contains~$\Lie
u_{\alpha+\beta}$, we must have $H = \{\, h \in N \mid h_{1,2} = p h_{2,3}
\, \}$, for some nonzero $p \in \real$. For each $t \in \real$, let $h_t =
h(1,1,1,pt,t,0)$. Then  $\rho(h_t) \asymp t^2 \asymp  \|h_t\|^2$, so $H$
is a Cartan-decomposition subgroup.

\begin{case}
 Assume that $\dim H \ge 3$.
 \end{case}
 We show that $H$ is a Cartan-decomposition subgroup.

If $\dim H \ge 4$, then $H$ must contain either $A$ or~$N$ \cf{HN=AN}, so
$H$ is a Cartan-decomposition subgroup. Thus, we may assume that $\dim H = 3$.
Furthermore, from Lemma~\ref{HN=AN}, we may assume that $\dim(H \cap N) = 2$.

 We may assume that $H \cap N$ is normalized by~$A$, for, otherwise,
Case~\ref{HinN} implies that $H \cap N$ (and, hence,~$H$) is a
Cartan-decomposition subgroup. So $\Lie h \cap \Lie n$ must be either
$\Lie u_\alpha + \Lie u_{\alpha+\beta}$ or $\Lie u_\beta + \Lie
u_{\alpha+\beta}$. For definiteness, let us assume that $\Lie h \cap \Lie n
= \Lie u_\alpha + \Lie u_{\alpha+\beta}$. (The calculations are
essentially the same in the other case. In fact, $\Lie u_\alpha + \Lie
u_{\alpha+\beta}$ and $\Lie u_\beta + \Lie u_{\alpha+\beta}$ are conjugate
under an outer automorphism of~$\Lie g$.)

\begin{subcase}
 Assume that $H = (H \cap A) \semiprod (H \cap N)$.
 \end{subcase}
 We may assume that $a \le 0$. (Recall that $a,b,c$ are defined
in~\eqref{SL3-abc}.) For each $t \in \real^+$, define $h^t = h(e^{at},
e^{bt}, e^{ct}, e^{dt},0,e^{dt}) \in H$, where $d = \max\{b,c\} > 0$. Then
$\rho(h^t) \asymp e^{2dt} \asymp \|h^t\|^2$, so $H$ is a
Cartan-decomposition subgroup.

\begin{subcase}
 Assume that $H \neq (H \cap A) \semiprod (H \cap N)$.
 \end{subcase}
 Let $\omega$ be a positive root as described in Lemma~\ref{not-semi}.
 Since $\Lie u_\alpha,\Lie u_{\alpha+\beta} \subset \Lie h$, we must have
$\omega = \beta$. Thus, for some nonzero~$p$, we have
 $$H = \{\, h(e^{t}, e^{t},e^{-2t}, u, pte^{t}, v) \mid t,u,v \in \real \,
\} .$$
 For each $t \in \real^+$, let 
 $ h_t = h(e^{t}, e^{t},e^{-2t}, pte^{t}, pte^{t}, 0) \in H$. Then
  $ \rho(h_t) \asymp t^2 e^{2t} \asymp \|h_t\|^2 $,
 so $H$ is a Cartan-decomposition subgroup.

\begin{case}
  Assume that $\dim H = 2$, that $H \not\subset N$, and that $H = (H \cap
A) \semiprod (H \cap N)$.
 \end{case}
 We may assume that $H \cap N$ is not a root group, for otherwise
Proposition~\ref{rootsemi-CDS} applies.
 We show that $H$ is a Cartan-decomposition subgroup if and only if $H
\cap A = \ker \alpha$ or $H \cap A = \ker \beta$. (That is, if and only
if $H$ is of type~\pref{SL3-CDS-alpha} or~\pref{SL3-CDS-beta},
respectively.) When $H$ is not a Cartan-decomposition subgroup, we have
$H \cap A = \ker(\alpha - \beta)$, in which case we show that $\rho(h)
\asymp h$ for every $h \in H$.

Because $H \cap N$ is normalized by $H \cap A$, and is one-dimensional,
every element of $\Lie h \cap \Lie n$ is an eigenvector for each element
of $H \cap A$. Because, by assumption, $\Lie h \cap \Lie n$ is not a root
space, this implies that two positive roots must agree on $H \cap A$. So $H
\cap A$ is the kernel of $\alpha$, $\beta$, or $\alpha - \beta$.

If  $H \cap A = \ker \alpha$ or $H \cap A = \ker \beta$, then we have
$\rho(a) \asymp \|a\|^2$ or $\rho(a^{-1}) \asymp \|a^{-1}\|^2$, for
every $a \in H \cap A$. Thus, $H$ is a Cartan-decomposition subgroup.

Assume now that $H \cap A = \ker(\alpha - \beta)$. Then
 $$ H = \{\, h(e^t,1,e^{-t},e^t u,pu, e^t p u^2/2) \mid t,u \in \real \,
\} ,$$
 for some nonzero $p \in \real$.
 Thus, we have $\rho(h) \asymp \max\{e^t u^2, e^{-t} \} \asymp h$ for every
$h \in H$ with $|t|,|u| >1$. Any element of~$H$ is within a bounded
distance of some point with $|t|,|u| > 1$, so we conclude from
Lemma~\ref{bdd-change} that $\mu(h) \asymp h$ for every $h \in H$.
Therefore, $H$ is not a Cartan-decomposition subgroup.

\begin{case}
 Assume that  $\dim H = 2$ and that $H \neq (H \cap A) \semiprod (H \cap
N)$.
 \end{case}
 We show that
 $H$ is not a Cartan-decomposition subgroup, and that 
 $$\mu(H) \asymp
\muH{\|h\|^{1/2} (\log \|h\|)}{\|h\|^2/(\log \|h\|)} .$$

Let the positive root~$\omega$ and the subgroup~$U$ of~$N$ be as described
in Lemma~\ref{not-semi}. Assume for definiteness that $\omega = \alpha$.
(The calculations are similar in the other cases. Indeed, the groups in
the other cases are conjugate to these under an automorphism of~$G$.)
Because $U_\alpha$ is one-dimensional, we must have $\psi(\ker \alpha) =
U_\alpha$.

Because the restrictions of $\beta$ and~$\alpha+\beta$ to $\ker\alpha$ are
nontrivial, unlike the restriction of~$\alpha$, and $\ker\alpha$
normalizes~$U$, but $U_\alpha \not\subset U$, we see that $U \subset
U_\beta U_{\alpha+\beta}$. Then, because $U = H \cap N$ is one-dimensional
and is normalized by~$U_\alpha$, we conclude that $H \cap N =
U_{\alpha+\beta}$. Thus, we have $H = \{
 h(e^{t},e^{t},e^{-2 t},pte^{t},0,s) \mid s,t \in\real
 \}$, for some nonzero~$p \in \real$. 
 For $|s|,|t| \ge 1$, we have
 $h \asymp \max\{ |t|e^t, e^{-2t},|s| \}$
 and $\rho \asymp \max\{ e^{2t}, |t|e^{-t}, |s|e^t \}$.

 Letting $s = t e^t \gg 0$ yields $\rho(h) \asymp t e^{2t} \asymp
\|h\|^2/(\log \|h\|)$. We now show that this is (approximately) the largest
possible size of $\rho(h)$ relative to~$h$. Because $\mu(H)$ is invariant
under the opposition involution \see{oppinv}, this implies that $\mu(H)
\asymp \muH{(\|h\| \log \|h\|)^{1/2}}{\|h\|^2/(\log \|h\|)}$, so $H$ is
not a Cartan-decomposition subgroup.
 For $t < -1$, we have $e^{2t} < 1$, $|t|e^{-t} < e^{-2t}$, and $|s|e^t <
|s|$, so $\|\rho(h)\| = O(\|h\|)$ is much smaller than $\|h\|^2/(\log
\|h\|)$. Now suppose that $t > 1$. If $|s| < |t|e^t$, then
 $$\rho(h) = O( t e^{2t}) = O \bigl( (t e^t)^2/t \bigr) =
O(\|h\|^2/\log\|h\|) .$$
 If $|s| \ge |t|e^t$, then
 $$\rho(h) = O( |s| e^t) = O \left( \frac{s^2}{|s|/e^t} \right)
 =
O \left( \frac{s^2}{\log s} \right)
 = O \left( \frac{\|h\|^2}{\log \|h\|} \right),$$
 as desired.
 \end{proof}

\begin{proof}[{\bf Proof of Theorem~\ref{SL3-minl}.}]
 Let $H$ be a minimal Cartan-decomposition subgroup of~$G =
\SL(3,\real)$. 

\setcounter{case}{0}

\begin{case} \label{SL3-minl-AN<3}
 Assume that $H \subset AN$ and that $\dim H < 3$.
 \end{case}
 Consider the possibilities given by Theorem~\ref{SL3-CDS}.
 If $H$ is of type~\fullref{SL3-CDS}{-A}, then $H = A$ is listed in
Theorem~\ref{SL3-minl}.
 If $H$ is of type~\fullref{SL3-CDS}{-inN}, then we may assume that $p =
1$, by replacing~$H$ with a conjugate via an element of~$A$; thus, $H$ is
the subgroup of~$N$ that is listed in Theorem~\ref{SL3-minl}.
 If $H$ is of type~\fullref{SL3-CDS}{-alpha}
or~\fullref{SL3-CDS}{-beta}, then we may assume that $p = 0$, by
replacing~$H$ with a conjugate via an element of~$\SO(2) \times \Id$
or~$\Id \times \SO(2)$, respectively (a maximal compact subgroup of the
centralizer of $H \cap A$); thus, $H$ is of
type~\fullref{SL3-CDS}{-root}, discussed below.

If $H$ is of type~\fullref{SL3-CDS}{-root}, then, because all roots are
conjugate under the Weyl group, we may assume that $\omega = \alpha$;
thus, $H$ is of the form~\pref{SL3-rootsemi}.
Example~\ref{SL3-rootsemi-eg} shows that we may assume that $\max\{p,q\}
= 1$ and that $\min\{p,q\} \ge -1/2$, so $H$ is listed in
Theorem~\ref{SL3-minl}.

\begin{case} \label{SL3-minl-AN=3}
 Assume that $H \subset AN$ and that $\dim H \ge 3$.
 \end{case}
 Because $H$ is minimal, we know that $H$ does not contain any conjugate
of~$A$, so $\dim(H \cap N) \ge 2$ \see{HN=AN}. Then, because $H \cap N$
is not a Cartan-decomposition subgroup, we see from
Corollary~\ref{SL3notCDS} that $\dim(H \cap N) = 2$ and that $\Lie h$
is either $\Lie u_\alpha + \Lie u_{\alpha+\beta}$ or $\Lie u_\beta +
\Lie u_{\alpha+\beta}$.  There is an automorphism~$\phi$ of~$G$ that
normalizes~$A$ and~$U_{\alpha+\beta}$, but interchanges~$U_\alpha$
with~$U_\beta$ (namely, $\phi$ is the composition of the Cartan
involution (``transpose-inverse") with the Weyl reflection corresponding
to the root~$\alpha+\beta$), so we may assume that $H = U_\beta +
U_{\alpha+\beta}$. 

If $H \cap A$ is nontrivial, then, because the Weyl chamber
containing~$A_\beta^+$ and the Weyl chamber
containing~$A_{\alpha+\beta}^+$ have no common interior, we see from
Proposition~\ref{rootsemi-CDS} that either $(H \cap A) \semiprod
U_\alpha$ or $(H \cap A) \semiprod U_{\alpha+\beta}$ is a
Cartan-decomposition subgroup, which contradicts the minimality of~$H$.

Thus, $H \cap A$ must be trivial, and we may assume that $H$ is
compatible with~$A$ \see{conj-to-compatible}, so, from
Lemma~\ref{not-semi}, we see that $H$ is conjugate (via an element
of~$A$) to the three-dimensional subgroup of~$AN$ that is listed in
Theorem~\ref{SL3-minl}. (We must have $\omega = \alpha$, because $H$
contains $U_\beta$ and~$U_{\alpha+\beta}$.)

\begin{case} \label{SL3-minl-notAN}
 Assume that $H$ is not conjugate to a subgroup of~$AN$. 
 \end{case}
 Any nontrivial, connected, semisimple Lie group has a connected,
proper, cocompact subgroup (namely, either the trivial subgroup or a
parabolic subgroup, depending on whether $H$ is compact or not), so the
minimality implies that $H$ is solvable.
 From the nilshadow construction (cf.~proof of Lemma~\ref{HcanbeAN}), we
may assume, after replacing~$H$ by a conjugate subgroup, that there is a
connected subgroup~$H'$ of~$AN$ and a compact, connected, abelian
subgroup~$T$ of~$G$, such that $T$ normalizes both~$H$ and~$H'$, and we
have $HT = H'T$. We may assume that $H'$ is compatible with~$A$
\see{conj-to-compatible}. Because $H' \cap N$ is a nontrivial, connected
subgroup of~$N$ that is normalized by a nontrivial, connected, compact
subgroup of~$G$, we must have either $H' \cap N = U_\alpha
U_{\alpha+\beta}$ or $H' \cap N = U_\beta U_{\alpha+\beta}$. Replacing
$H$ by a conjugate via an automorphism of~$G$ (cf.~the
automorphism~$\phi$ described in Case~\ref{SL3-minl-AN=3}), we may assume
that $H' \cap N = U_\beta U_{\alpha+\beta}$. Then $H' \cap N$ is not a
Cartan-decomposition subgroup, so $H' \not\subset N$. Let $T_H = (HN)
\cap A$. From Lemma~\ref{not-semi}, we know that $T_H$ is contained in
the Zariski closure of~$H'$, so, by replacing~$H$ with a conjugate, we
may assume that $T$ centralizes~$T_H$. Thus, $T_H = \ker \omega$, for
some root~$\omega$, and $T \subset \langle U_\omega, U_{-\omega}\rangle$.
Because $U_{\beta} U_{\alpha+\beta} = H' \cap N$ is normalized by~$T$,
we must have $\omega = \alpha$.

Now $T T_H$ is a maximal connected subgroup of the centralizer
$C_G(T_H)$ (because $\SO(2)$ is maximal in $\SL(2,\real)$), so 
 $$ H' \cap C_G(T_H) = (H' T T_H) \cap C_G(T_H) \cap H' = (T T_H)
\cap H' \subset T_H .$$
 Thus, $H'$ is not of type~\fullref{not-semi}{-not}, so
Lemma~\ref{not-semi} implies that $H' = (H' \cap A) \semiprod (H' \cap
N)$. Therefore, $H \subset (T \ker \alpha) \semiprod (U_\beta
U_{\alpha+\beta})$, so, after conjugating $T$ to a subgroup of $\SO(2)
\times \Id$, it is easy to see that $H$ is of the form~\pref{SL3-rot}.
 \end{proof}

\section{Subgroups of $\SO(2,n)$ contained in~$N$}
 \label{SO2n-HinN-sect}

We now study $G = \SO(2,n)$. In this section, we determine which subgroups
of~$N$ are Cartan-decomposition subgroups of $G$.
Theorem~\ref{SO2n-HinN-CDS} gives simple conditions to check whether a
subgroup of~$N$ is a Cartan-decomposition subgroup, but the statement of
the result requires some notation. Theorem~\ref{HinN-notCDS} describes
all the subgroups of~$N$ that are \emph{not} Cartan-decomposition
subgroups, and Proposition~\ref{HinN-Cartanproj} describes the image of
each of these subgroups under the Cartan projection. Analogous results
for subgroups not contained in~$N$ appear in
Section~\ref{SO2n-notinN-sect}.

\begin{notation}
 Recall that we realize $\SO(2,n)$ as isometries of the indefinite form
$\langle v \mid v \rangle = v_1 v_{n+2} + v_2 v_{n+1} + \sum_{i=3}^n
v_i^2$ on~$\real^{n+2}$, that $A$ consists of the group of diagonal
matrices, and that $N$ consists of the upper-triangular unipotent
matrices \see{SO2n-defn}.
 Let 
 $$A^+= \{\, \diag (a_1, a_2 , 1,1,\ldots,1,1 , a_2^{-1} , a_1^{-1})\mid 
a_1 \geq a_2 \geq 1 \,\} $$
 and
 $K = \SO(2,n) \cap \SO(n+2)$.
 \end{notation}

\begin{notation} \label{h(phi,x,y,eta)}
 Given $\phi,\eta \in \real$ and $x,y \in \real^{n-2}$ (and letting $t_1
= t_2 = 0$), there is a corresponding element of~$\Lie n$ \see{SO2n-AN}.
The exponential of this is an element~$h = h(\phi,x,y,\eta)$ of~$N$.
Namely,
 $$ h = \begin{pmatrix}
 1 & \phi & x +\frac{1}{2} \phi y & \eta -\frac{1}{2} (x \cdot
y) -\frac 16 \phi \|y\|^2 & -\phi \eta -\frac{1}{2} \|x\|^2 +
\frac{1}{24} \phi^2 \|y\|^2 
\\
 &1& y &-\frac{1}{2} \|y\|^2 & -\eta
 -\frac{1}{2} (x \cdot y) +\frac 16 \phi \|y\|^2\\ 
 &  & \Id & -y^\transpose & -x^\transpose  +\frac{1}{2} \phi y^\transpose \\
 && &1& -\phi \\
 &&& &1 
 \end{pmatrix} .$$
 Because the exponential map is a diffeomorphism from~$\Lie n$ onto~$N$,
each element of~$N$ has a unique representation in this form. Thus, each
element of~$N$ determines corresponding values of~$\phi$, $x$, $y$,
and~$\eta$. We sometimes write $\phi_h$, $x_h$, $y_h$,
and~$\eta_h$ for these values, to emphasize the element~$h$ of~$N$ that is
under consideration.
 \end{notation}

The following description of the Cartan-decomposition subgroups contained
in~$N$ is obtained by combining Lemmas~\ref{HinN-linear-dim>1}
and~\ref{HinN-square} with Proposition~\ref{dim>Rrank}.

\begin{thm} \label{SO2n-HinN-CDS}
 Assume that $G = \SO(2,n)$. A closed, connected subgroup~$H$ of~$N$ is a
Cartan-decomposition subgroup of~$G$ if and only if either
 \begin{enumerate}
 \item \label{SO2n-HinN-CDS-dim2}
 $\dim H = 2$, $\Lie u_{\alpha+2\beta} \subset \Lie h$, and there
is an element of~$\Lie h$ such that $\phi \neq 0$ and $y \neq 0$; or
 \item \label{SO2n-HinN-CDS-dim>2}
 $\dim H \ge 2$ and
 \begin{enumerate}
 \item there is an element of~$\Lie h$ such that $\langle (\phi,x),
(0,y)\rangle$ is one-dimensional in $\real^{n-1}$; and
 \item
 there is a nonzero element of~$\Lie h$ such that either 
 \begin{enumerate}
 \item \label{SO2n-HinN-CDS-<xy>=2}
 $\langle (\phi,x), (0,y)\rangle$ is two-dimensional in $\real^{n-1}$; or 
 \item \label{SO2n-HinN-CDS-y=0}
 $y = 0$ and $\|x\|^2 = -2 \phi \eta$.
 \end{enumerate}
 \end{enumerate}
 \end{enumerate}
 \end{thm}

\begin{notation}
 We let $\alpha$ and~$\beta$ be the simple real roots of $\SO(2,n)$,
defined by $\alpha(a) = a_1/a_2$ and $\beta(a) = a_2$, for an element~$a$
of~$A$ of the form
 $$ a= \diag (a_1, a_2 , 1,1,\ldots,1,1 , a_2^{-1} ,
a_1^{-1}) . $$
 Thus, 
 \begin{itemize}
 \item the root space $\Lie u_\alpha$ is the $\phi$-axis in~$\Lie n$, 
 \item the root space $\Lie u_\beta$ is the $y$-subspace in~$\Lie n$, 
 \item the root space $\Lie u_{\alpha+\beta}$ is the $x$-subspace in~$\Lie n$, and
 \item  the root space $\Lie u_{\alpha+2\beta}$ is the $\eta$-axis
in~$\Lie n$.
 \end{itemize}
 \end{notation}

Theorem~\ref{SO2n-HinN-CDS} describes Cartan-decomposition subgroups. The
following result describes the subgroups of~$N$ that are \emph{not}
Cartan-decomposition subgroups. It is obtained by combining Corollaries
\ref{HinN-nolinear} and~\ref{HinN-nosquare} with
Proposition~\ref{dim>Rrank}.

\begin{thm} \label{HinN-notCDS}
 Assume that $G = \SO(2,n)$. A closed, connected subgroup~$H$ of~$N$ is not a
Cartan-decomposition subgroup of~$G$ if and only if either
 \begin{enumerate}

\item \label{HinN-dim1}
 $\dim H \le 1$; or

 \item \label{HinN-<xy>not1}
 for every nonzero element of~$\Lie h$, we have $\phi = 0$ and $\dim
\langle x,y \rangle \neq 1$ {\upshape(}i.e., $\dim \langle x,y \rangle \in
\{0,2\}${\upshape)}; or

 \item \label{HinN-<xy>=1}
 for every nonzero element of~$\Lie h$, we have $\phi = 0$ and $\dim
\langle x,y \rangle = 1$; or

 \item \label{HinN-b2-2a}
 there exists a subspace~$X_0$ of~$\real^{n-2}$, $b \in X_0$, $c \in
X_0^\perp$, and $p \in \real$ with $\|b\|^2 - \|c\|^2 - 2p < 0$, such that
for every element of~$\Lie h$, we have $y= 0$, $x \in \phi c + X_0$, and
$\eta = p\phi + b \cdot x$ {\upshape(}where $b \cdot x$ denotes the
Euclidean dot product of the vectors $b$ and~$x$
in~$\real^{n-2}${\upshape)}.
 \end{enumerate}
 \end{thm}

In the course of the proof of Theorem~\ref{HinN-notCDS}, we
calculate the Cartan projection of  each subgroup that is not a
Cartan-decomposition subgroup. Proposition~\ref{HinN-Cartanproj}
collects these results. The statement (and all of the proofs in this
section) is based on Proposition~\ref{CDS<>seq}, so we describe the
required representations $\rho_1$ and~$\rho_2$ of~$G$.

\begin{notation} \label{SO2n-rho}
 Define a representation 
 $$\rho \colon \SO(2,n) \to \SL(\real^{n+2} \wedge \real^{n+2})
 \mbox{ by }
 \rho(g) = g \wedge g .$$
 \end{notation}
 In the notation of Proposition~\ref{norm}, we have
 \begin{eqnarray*}
 \rho _1 &=& \mbox {the standard representation on~$\real^{n+2}$
 \qquad
($\rho_1(g) = g$)}, \\
 \rho _2 &=& \rho \qquad (= \rho_1 \wedge \rho_1), \\
 k_1 &=& 1 \mbox{ \qquad and \qquad } k_2=2 .
 \end{eqnarray*}
 Thus, Proposition~\ref{norm} yields the following
fundamental lemma.

\begin{lem}
 Assume that $G = \SO(2,n)$.
 A subgroup~$H$ of~$G$ is a Cartan-decomposition subgroup if and only if
$\mu(H) \approx \muH{\|h\|}{\|h\|^2}$, where the representation~$\rho$ is
defined in  Notation~\ref{SO2n-rho}.
 \end{lem}

 \begin{prop} \label{HinN-Cartanproj}
 Assume that $G = \SO(2,n)$ and let $H$ be a closed, connected subgroup of~$N$.
 \begin{enumerate}
 \item 
 If $H$ is of type~\ref{HinN-notCDS}\pref{HinN-dim1} {\upshape(}i.e., if
$\dim H \le 1${\upshape)}, then either $\rho(h) \asymp h$ for all $h \in
H$, or $\rho(h) \asymp \|h\|^{3/2}$ for all $h \in
H$, or $\rho(h) \asymp \|h\|^2$ for all $h \in H$.
 \item 
 If $H$ is of type~\ref{HinN-notCDS}\pref{HinN-<xy>not1}, then $\rho(h)
\asymp  \|h\|^2$ for all $h \in H$.
 \item 
 If $H$ is of type~\ref{HinN-notCDS}\pref{HinN-<xy>=1}, then $\rho(h)
\asymp h$ for all $h \in H$.
 \item 
 If $H$ is of type~\ref{HinN-notCDS}\pref{HinN-b2-2a}, then $\rho(h)
\asymp h$ for all $h \in H$.
 \end{enumerate}
 \end{prop}

Note that if either  $\rho(h) \asymp h$ for all $h \in H$, or $\rho(h)
\asymp \|h\|^2$ for all $h \in H$, then $\mu(H)$ is within a bounded
distance of one of the walls of~$A^+$. On the other hand, if $H$ is a
Cartan-decomposition subgroup, then $\mu(H)\approx A^+$. Thus, certain
one-dimensional subgroups (described in Lemma~\ref{rho(phiy)}) are the
only connected subgroups of~$N$ for which $\mu(H)$ is neither a wall of
the Weyl chamber nor all of~$A^+$ (up to bounded distance). In
Section~\ref{SO2n-notinN-sect}, where $H$ is not assumed to be contained
in~$N$, we will see several more examples of this.

\begin{proof}[{\bf Proof of Theorem~\ref{SO2n-inN-minl}.}]
 It is easy to see from Theorem~\ref{SO2n-HinN-CDS} that the subgroup
of~$G$ corresponding to each of the subalgebras listed in
Theorem~\ref{SO2n-inN-minl} is a Cartan-decomposition subgroup. Each is
minimal, because each is 2-dimensional \see{dim>Rrank}.
Lemma~\ref{Nnotconj} shows that no two of the listed subalgebras are
conjugate.

Thus, given a minimal
Cartan-decomposition subgroup of~$G$ that is contained in~$N$, all that
remains is to show that the Lie algebra of~$H$ is conjugate to one on
the list.

 The Weyl reflection corresponding to the root~$\beta$
fixes~$\alpha+\beta$, but interchanges~$\alpha$ and~$\alpha+2\beta$.
Therefore, the subalgebra
  \begin{equation} \label{SO2n-inN-minl-a+b,a+2b}
  \bigset{
 \begin{pmatrix}
 0&0 &x&0&0 &\eta&0 \\
  &0 &0&0&0 &0&-\eta \\
  &&& \cdots \\
 \end{pmatrix}
 }{ x,\eta \in \real}
 \end{equation}
 is conjugate to the subalgebra of type~\fullref{SO2n-inN-minl}{-a,a+b}
with $\epsilon_2 = 0$.

\setcounter{case}{0}

 \begin{case}
 Assume that $\Lie h \subset \Lie u_\beta + \Lie u_{\alpha+\beta} + \Lie
u_{\alpha+2\beta}$.
 \end{case}
 From Theorem~\ref{SO2n-HinN-CDS}, we know that there is an element~$u$
of~$\Lie h$ with $\dim\langle x_u,y_u\rangle = 1$. By replacing~$H$ with
a conjugate under~$\langle U_\alpha, U_{-\alpha}\rangle$, we may assume
that $y_u = 0$. Then, by replacing~$H$ with a conjugate under~$U_\beta$,
we may assume that $\eta_u = 0$. If $\Lie u_{\alpha+2\beta} \subset \Lie
h$, then $\Lie h$ contains a conjugate of the
subalgebra~\pref{SO2n-inN-minl-a+b,a+2b}. Thus, we may now assume
that $\Lie u_{\alpha+2\beta} \not\subset \Lie h$. Therefore, $\Lie h$
has no nonzero elements of type~\fullref{SO2n-HinN-CDS}{-y=0}, so there
must be some $v \in H$ with $\dim \langle x_v,y_v \rangle = 2$. Because
$\Lie u_{\alpha+2\beta} \not\subset \Lie h$, we know that $\Lie h$ is
abelian, so $v$ commutes with~$u$, which means that $x_u$ is
perpendicular to~$y_v$. Thus, $\langle u,v \rangle$ is conjugate to one
of the two subalgebras of type~\fullref{SO2n-inN-minl}{-a+b,b}. 

\begin{case}
 Assume that there is an element~$u$ of~$\Lie h$, such that $u \notin
\Lie u_\beta + \Lie u_{\alpha+\beta} + \Lie u_{\alpha+2\beta}$.
 \end{case}

\begin{subcase} \label{SO2n-minlpf-b}
 Assume that there exists
 $v \in \bigl( \Lie h \cap (\Lie u_\beta + \Lie u_{\alpha+\beta} + \Lie
u_{\alpha+2\beta}) \bigr) \setminus (\Lie u_{\alpha+\beta} + \Lie
u_{\alpha+2\beta})$.
 \end{subcase}
 We have $[u,v] \in (\Lie u_{\alpha+\beta} + \Lie
u_{\alpha+2\beta}) \setminus \Lie u_{\alpha+2\beta}$ and
 $[u,v,v] \in \Lie u_{\alpha+2\beta} \setminus \{0\}$, 
 so $\langle [u,v], [u,v,v] \rangle$ is conjugate to the
subalgebra~\pref{SO2n-inN-minl-a+b,a+2b}.

\begin{subcase} \label{SO2n-minlpf-a,a+b,a+2b}
 Assume that $\Lie h \subset \Lie u_\alpha + \Lie u_{\alpha+\beta} + \Lie
u_{\alpha+2\beta}$.
 \end{subcase} 
 No element of~$\Lie h$ is of type~\fullref{SO2n-HinN-CDS}{-<xy>=2}, so
there must be a nonzero element~$v$ of~$\Lie h$ with $\|x_v\|^2 =
-2\phi_v\eta_v$. Replacing $H$ by a conjugate under~$U_\beta$, we may
assume that $x_v = 0$; then either $\phi_v$ or~$\eta_v$ must also
be~$0$. (Even after conjugation, we must have $\rho \bigl( \exp (tv)
\bigr) \asymp \|\exp (tv)\|^2$.)

 If $\phi_v = 0$, then $\Lie u_{\alpha+2\beta} \subset \Lie h$, so $\Lie
h$ contains a conjugate of either the
subalgebra~\pref{SO2n-inN-minl-a+b,a+2b} or the
subalgebra~\fullref{SO2n-inN-minl}{-a,a+2b}.

 If $\eta_v = 0$, then $\Lie u_{\alpha} \subset \Lie h$. Let $w$ be a
nonzero element of $\Lie h \cap (\Lie u_{\alpha+\beta} + \Lie
u_{\alpha+2\beta})$. If $\eta_w = 0$, then $\langle v,w \rangle$ is
conjugate to the subalgebra of type~\fullref{SO2n-inN-minl}{-a,a+b} with
$\epsilon_2 = 0$. If $\eta_2 \neq 0$, then, by replacing $\Lie h$ with a
conjugate under~$U_{-\beta}$, we may assume that $x_w = 0$, so 
$\langle v,w \rangle$ is the subalgebra~\fullref{SO2n-inN-minl}{-a,a+2b}.

\begin{subcase}
 The general case.
 \end{subcase}
 If $H$ is of type~\fullref{SO2n-HinN-CDS}{-dim2}, then $\Lie h$ is
conjugate to the subalgebra of type~\fullref{SO2n-inN-minl}{-a,a+2b} with
$\epsilon_1 = 1$. Thus, we henceforth assume that $H$ is of
type~\fullref{SO2n-HinN-CDS}{-dim>2}.

 Let $v$ be an element of~$\Lie h$ with $\dim \langle (\phi_v,x_v),
(0,y_v) \rangle = 1$. We may assume that
Subcase~\ref{SO2n-minlpf-a,a+b,a+2b} does not apply, so we may assume
that $y_u \neq 0$.

Suppose that $\phi_v \neq 0$. Then $y_v$ must be~$0$. Because $\phi_v \neq 0$, $y_v = 0$ and $y_u
\neq 0$, there is some linear combination of~$v$ and~$u$, such that
$\phi = 0$ and $y \neq 0$, so Subcase~\ref{SO2n-minlpf-b} applies.

We may now assume that $\phi_v = 0$. We may assume that $y_v$ is
also~$0$, for, otherwise, Subcase~\ref{SO2n-minlpf-b} applies. Thus, $v
\in (\Lie u_{\alpha+\beta} + \Lie u_{\alpha+2\beta}) \setminus \Lie
u_{\alpha+2\beta}$. We may assume that $[u,v] =0$; otherwise,
$\langle v, [u,v] \rangle$ is conjugate to the
subalgebra~\pref{SO2n-inN-minl-a+b,a+2b}, because $[u,v] \in
\Lie u_{\alpha+2\beta}$.

Conjugating by an element of~$U_{\beta}$, we may assume that $x_u$ is
a scalar multiple of~$x_v$, and that $\eta_v = 0$. Then, by
replacing~$u$ with $u + \lambda v$, for an appropriate scalar~$\lambda$,
we may assume that $x_u = 0$. Recall that $y_u \neq 0$. Therefore,
conjugating by an element of~$U_{\alpha+\beta}$, we may assume that
$\eta_u = 0$. Because $[u,v] = 0$, we know that $x_v$
is perpendicular to~$y_u$. Thus, $\langle u,v \rangle$ is conjugate to
one of the two subalgebras of type~\fullref{SO2n-inN-minl}{-a,a+b}.
 \end{proof}

\begin{lem} \label{Nnotconj}
 Assume that $G = \SO(2,n)$.
 No two of the subalgebras of~$\Lie n$ listed in
Theorem~\ref{SO2n-inN-minl} are conjugate under~$G$. {\upshape(}In fact,
they are not even conjugate under $\GL(n+2,\real)$.
 \end{lem}

\begin{proof} For each subalgebra~$\Lie h$, we look at the
restriction of the exponential map to~$\Lie h$. This is a polynomial
function from~$\Lie h$ into~$N$. For convenience, we use $(\phi,\eta)$,
$(x,y)$, or $(\phi,x)$ as coordinates on~$\Lie h$, and use the matrix
entries as coordinates on~$N$.

The subalgebra of type~\fullref{SO2n-inN-minl}{-a,a+2b} with $\epsilon_1
= 1$ and the subalgebra of type~\fullref{SO2n-inN-minl}{-a,a+b} with
$\epsilon_2 = 1$ are the only cases where the exponential has a 4th
degree term (namely, $\phi^4/24$), so they cannot be conjugate to any of
the others. They are not conjugate to each other, because in the
subalgebra of type~\fullref{SO2n-inN-minl}{-a,a+2b} with $\epsilon_1 =
1$, there is a nontrivial subspace~$\Lie u$ such that $\mu(\exp \Lie u)
\approx \muH{\|h\|^2}{\|h\|^2}$ (namely, $\Lie u = \Lie
u_{\alpha+2\beta}$), but there is no such subspace in the subalgebra of
type~\fullref{SO2n-inN-minl}{-a,a+b} with $\epsilon_2 = 1$.

Of the others, it is only for the two subalgebras of
type~\fullref{SO2n-inN-minl}{-a+b,b} that there are two independent
quadratic terms (so the polynomial $\rho\bigl( \exp(u) \bigr)$ is of
degree~4).  For the subalgebra with $\epsilon_3 = 0$, one or the other of
the quadratic terms vanishes if either~$x$ or~$y$ is~$0$ (that is, on
two lines); whereas, for $\epsilon_3 = 1$, neither quadratic term vanishes
unless $y = 0$ (a single line).

For each of two remaining subalgebras, the exponential function has a
single quadratic term. For the subalgebra of
type~\fullref{SO2n-inN-minl}{-a,a+b} with $\epsilon_2 = 0$,
the quadratic term does not vanish unless $x = 0$ (a single line).
Whereas, for the subalgebra of type~\fullref{SO2n-inN-minl}{-a,a+2b}
with $\epsilon_1 = 0$, the quadratic term vanishes if either $\phi= 0$
or~$\eta = 0$ (two lines).
 \end{proof}

The remainder of this section has three parts. In
\S\ref{HinN-linear-sect}, we determine when $\mu(H)$ contains elements
close to the wall given by $\rho(h) \asymp h$. In
\S\ref{HinN-square-sect}, we determine when $\mu(H)$ contains elements
close to the other wall, given by $\rho(h) \asymp \|h\|^2$. In
\S\ref{HinN-calculate}, we calculate $\mu(H)$ for some subgroups that are
not Cartan-decomposition subgroups.

\subsection{When is the size of $\rho(h)$ linear?}
\label{HinN-linear-sect}

\begin{lem} \label{HinN-linear-dim>1}
 Assume that $G = \SO(2,n)$.
 Suppose that $H$ is a closed, connected subgroup of~$N$, and that $\dim
H \neq 1$. Then the following are equivalent:
 \begin{enumerate}
 \item There is a sequence $\{h_n\}$ of elements of~$H$ with $h_n \to
\infty$, such that $\rho(h_n) \asymp h_n$.
 \item There is a continuous curve $\{h_t\}_{t \in \real^+}$ of elements
of~$H$ with $h_t \to \infty$, such that $\rho(h_t) \asymp h_t$.
 \item Either
 \begin{enumerate}
 \item there is an element of~$\Lie h$ such that $\langle (\phi,x),
(0,y)\rangle$ is one-dimensional in $\real^{n-1}$; or
 \item  $\dim H = 2$, $\Lie u_{\alpha+2\beta} \subset \Lie h$, and there
is an element of~$\Lie h$ such that $\phi \neq 0$ and $y \neq 0$.
 \end{enumerate}
 \end{enumerate}
 \end{lem}

We prove the following more general (but slightly more complicated)
version that does not assume that $\dim H \neq 1$.

\begin{lem} \label{HinN-linear}
  Assume that $G = \SO(2,n)$, and let $H$ be a closed, connected subgroup
of~$N$.
Then the following are equivalent:
 \begin{enumerate}
 \item \label{HinN-linear-seq}
 There is a sequence $\{h_m\}$ of elements of~$H$ with $h_m \to
\infty$, such that $\rho(h_m) \asymp h_m$.
 \item \label{HinN-linear-curve}
 There is a continuous curve $\{h_t\}_{t \in
\real^+}$ of elements of~$H$ with $h_t \to \infty$, such that $\rho(h_t)
\asymp h_t$.
 \item \label{HinN-linear-Liealg}
 Either
 \begin{enumerate}
 \item \label{HinN-linear-Liealg-1D}
 There is an element of~$\Lie h$ such that $\langle (\phi,x),
(0,y)\rangle$ is one-dimensional in $\real^{n-1}$, and either $y \neq 0$
or $\|x\|^2 \neq -2\phi\eta$; or
 \item  \label{HinN-linear-Liealg-dim2}
 $\dim H = 2$, $\Lie u_{\alpha+2\beta} \subset \Lie h$, and there
is an element of~$\Lie h$ such that $\phi \neq 0$ and $y \neq 0$.
 \end{enumerate}
 \end{enumerate}
 \end{lem}

\begin{proof}
 (\ref{HinN-linear-seq} $\Rightarrow$ \ref{HinN-linear-Liealg})

\setcounter{case}{0}

\begin{case}
 Assume that $\phi = O(1)$ in $h_m$.
 \end{case}
 We may assume that $\phi = 0$ \see{bdd-change}.
  Let $a_m =
k'_m h_m k_m$ be the Cartan projection of~$h_m$, with $k'_m,k_m \in K$, and
write $a_m = \diag(a_m^{(1)},\ldots,a_m^{(n+2)})$, with $a_m^{(1)} \ge
a_m^{(2)} \ge \cdots \ge a_m^{(n+2)}$ (because $a_m \in A^+$). Let 
 $$ V = \{\, v \in \real^{n+2} \mid \forall i \le n, ~ v_i = 0 \, \}
 \mbox{ and }
 W = \{\, w \in \real^{n+2} \mid w_1 = 0 \, \}
 ,$$
 so $\dim V = 2$ and $\dim W = n+1$.
 Then $\dim V + \dim (k_m W) = n+3 > n+2$, so there is a nonzero vector
$v_m \in V \cap (k_m W)$; we may assume that $\|v_m\| = 1$.

Because $v_m \in V$, we may write $v_m = (0,0,\ldots,0,q_m,p_m)$, so 
 $$h_m v_m = \bigl(*,*, -(p_m x_m + q_m y_m),*,* \bigr) .$$
 Because $\rho(a_m) \asymp \rho(h_m) \asymp h_m \asymp a_m$, we must have
$a_m^{(i)} = O(1)$ for $i \ge 2$, so $a_m w = O(w)$ for all $w \in W$.
Then, since $k_m^{-1} v_m \in W$, we have
 \begin{eqnarray*}
 \| p_m x_m + q_m y_m \|
 &\le& \| h_m v_m \|
 = \| (k'_m)^{-1} a_m k_m^{-1} v_m \| \\
 &=& \| a_m k_m^{-1} v_m \|
 = O(k_m^{-1} v_m)
 = O(v_m)
 = O(1) .
 \end{eqnarray*}
 Passing to a subsequence, we may assume that $\{v_m\}$
converges, so $\{p_m\}$~and~$\{q_m\}$ converge. Clearly, $\|x_m\| +
\|y_m\|$ tends to infinity, because $\rho(h) \asymp \|h\|^2$ if $\phi$,
$x$, and~$y$ are all zero (so the same is true if $\phi,x,y = O(1)$
\see{bdd-change}). Thus, letting $x'_m = x_m/(\|x_m\| + \|y_m\|)$ and $y'_m
= y_m/(\|x_m\| + \|y_m\|)$, we have $p_m x'_m + q_my'_m = o(1)$. Because
$\|x'_m\|,\|y'_m\|\le 1$, we may pass to a convergent subsequence. Then
$p_\infty x'_\infty + q_\infty y'_\infty = 0$, so $x'_\infty$
and~$y'_\infty$ are linearly dependent. Because $\Lie h$ is a closed
subset of~$\Lie n$, there is some element of~$\Lie h$ with $\phi = 0$, $x
= x'_\infty$, and $y = y'_\infty$, so
conclusion~\pref{HinN-linear-Liealg-1D} holds.

\begin{case}
 Assume that $\phi$ is unbounded in~$\{h_m\}$.
 \end{case}
 By passing to a subsequence, we may assume that $|\phi_m| \to \infty$.

Assume for the moment that some element~$u$ of $\Lie h \cap (\Lie u_\beta
\oplus \Lie u_{\alpha +\beta } \oplus \Lie u_{\alpha +2\beta })$ does not
belong to $\Lie u_{\alpha +2\beta }$. If $u \notin \Lie u_{\alpha +\beta }
\oplus \Lie u_{\alpha +2\beta }$, then $[u, \log h]$ belongs to $\Lie
u_{\alpha +\beta } \oplus \Lie u_{\alpha +2\beta }$, but is not in $\Lie
u_{\alpha +2\beta }$. Thus, by replacing $u$ with $[u, \log h]$ if
necessary, we may assume, without loss of generality, that $u \in \Lie
u_{\alpha +\beta } \oplus \Lie u_{\alpha +2\beta }$. This implies that
$\phi_u = 0$, $y_u = 0$ and $x_u \neq 0$, so $\dim \langle (\phi_u,x_u),
(0,y_u) \rangle = 1$, as desired.

We may henceforth assume that every element of $\Lie h \cap (\Lie u_\beta \oplus
\Lie u_{\alpha +\beta } \oplus \Lie u_{\alpha +2\beta })$ belongs to $\Lie
u_{\alpha +2\beta }$. This means that $x$~and~$y$ are uniquely determined
by~$\phi$, so $\dim H \le 2$,  and there are ${\hat \phi} \in \Lie u_\alpha $,
${\hat x} \in \Lie u_{\alpha +\beta }$, and ${\hat y} \in \Lie u_{\beta }$, such
that ${\hat \phi} \neq 0$, and $\Lie h \subset \real({\hat \phi} + {\hat x} + {\hat y}) + \Lie
u_{\alpha +2\beta }$. 

Suppose that $\hat y \neq 0$. If $\dim H = 1$, then
Lemma~\ref{rho(phiy)} implies that $\rho(h) \asymp \|h\|^{3/2}$  for $h
\in H$, so $\rho(h) \not\asymp h$, which is a contradiction. Thus, we
must have $\dim H = 2$, so conclusion~\pref{HinN-linear-Liealg-dim2}
holds.

We now know ${\hat y} = 0$. If $\dim \Lie h = 1$ and $\|x\|^2 = -2 \phi
\eta$ for some (and, hence, all) elements of~$\Lie h$, then
Lemma~\fullref{HinN-square}{-x^2=2PhiEta} implies that $\rho(h)
\asymp \|h\|^2 \not\asymp h$, a contradiction. If $\dim \Lie h = 2$,
then $\Lie u_{\alpha +2\beta } \subset \Lie h$, so it is easy to choose
an element of~$\Lie h$ with $\|x\|^2 \neq -2 \phi \eta$.

\medbreak

 (\ref{HinN-linear-Liealg} $\Rightarrow$ \ref{HinN-linear-curve})

\setcounter{case}{0}

\begin{case}
 Assume that there is an element~$h$ of~$H$, such that $\dim \langle
(\phi,x), (0,y) \rangle = 1$, and $y \neq 0$.
 \end{case}
 Because $\dim \langle (\phi,x), (0,y) \rangle = 1$, we know that
$(\phi,x)$ is a scalar multiple of $(0,y)$, so $\phi =0$. Furthermore,
we may assume, after replacing~$H$ (and, hence,~$h$) with a conjugate
under~$U_{\alpha}$,
that $x = 0$ \see{conj-alm-same}. Then, letting $h^t = \exp(t \log h)$,
we have 
 $$ h^t_{ij} =
 \begin{cases}
 O(1) & \mbox{if $i \neq 2$ and $j \neq n+1$} \\
 O(t) & \mbox{if $(i,j) \neq (2,n+1)$} \\
 \end{cases}
 $$
 and $h^t_{2 , n+1} \asymp t^2$. Therefore $\rho(h^t) \asymp t^2 \asymp
h^t$, as desired.

\begin{case} \label{HinNcurve-x^2=-2phieta}
  Assume that there is an element~$h$ of~$H$, such that $y = 0$ and
$\|x\|^2 \neq -2\phi\eta$.
 \end{case}
 Letting $h^t = \exp(t \log h)$, we have 
 $$ h^t_{ij} =
 \begin{cases}
 O(1) & \mbox{if $i \neq 1$ and $j \neq n+2$} \\
 O(t) & \mbox{if $(i,j) \neq (1,n+2)$} \\
 \end{cases}
 $$
 and, because $\|x\|^2 \neq -2 \phi \eta$, we have $h^t_{1 , n+2} \asymp
t^2$. Therefore it is not difficult to see that $\rho(h^t) \asymp t^2
\asymp h^t$, as desired.

\begin{case}
 Assume that $\dim H = 2$, $\Lie u_{\alpha+2\beta} \subset \Lie h$,
and there is an element~$u$ of~$\Lie h$ such that $\phi_u \neq 0$
and $y_u \neq 0$.
 \end{case}
 Replacing $H$ by a conjugate under~$U_\beta$, we may assume that
$x_u = 0$. For any large real number~$t$, let $h = h^t$ be the
element of $\exp(t u + \Lie u_{\alpha+2\beta})$ that satisfies $\eta =
-\phi \|y\|^2/{12}$. Then
 $$ h = 
 \begin{pmatrix}
 1 & \phi & \frac{1}{2} \phi y & - \frac{1}{4} \phi \|y\|^2
   & \frac{1}{8} \phi^2 \|y\|^2 \\
   &  1   & y                  & - \frac{1}{2} \|y\|^2
   &  \frac{1}{4} \phi \|y\|^2 \\
   &      &  1                 & -y^\transpose      
   &  \frac{1}{2} \phi y^\transpose   \\
   &      &                    &  1                   &   -\phi \\
   &      &                    &                      &  1 \\
 \end{pmatrix}
 .$$
 Because $t$ is large, we know that $\phi$ and~$\|y\|$ are large,
so it is clear that $h \asymp \phi^2 \|y\|^2$.

Let $h'$ be the matrix obtained from~$h$ by deleting the
first two columns. Then the first two rows
of~$h'$ are linearly dependent (because the first row of~$h'$ is
$\phi/2$ times the remainder of the second row). Therefore, we have 
 $$\det \begin{pmatrix}
 h_{1,i} & h_{1,j} \\ h_{2,i} & h_{2,j} \\
 \end{pmatrix}
 = 0,
  \qquad \mbox{whenever $i,j > 2$.}
 $$
 Similarly,  we have 
 $$\det \begin{pmatrix}
 h_{i,n+1} & h_{i,n+2} \\ h_{j,n+1} & h_{j,n+2} \\
 \end{pmatrix}
 = 0,
  \qquad \mbox{whenever $i,j \le n$.}
 $$
 It is easy to see that the determinant of any
other $2 \times 2$ submatrix of~$h$ is $O(\phi^2 \|y\|^2)$. Thus, we
conclude that $\rho(h) = O(\phi^2 \|y\|^2) = O(h)$, as desired.
 \end{proof}

\begin{cor}[of proof] \label{HinN-nolinear}
 Assume that $G = \SO(2,n)$, and let $H$ be a closed, connected subgroup of~$N$.
The following are equivalent:
 \begin{enumerate}
 \item
 There is no sequence $\{h_m\}$ of elements of~$H$ with $h_m \to \infty$
and $\rho(h_m) \asymp h_m$.
 \item
 There is no curve $\{h_t\}_{t \in \real^+}$ in~$H$ with $h_t
\to \infty$ and $\rho(h_t) \asymp h_t$.
 \item
 Either
 \begin{enumerate}
 \item \label{phi=0,rank1}
 for every element of~$\Lie h$, we have $\phi = 0$ and $\dim \langle
x,y \rangle \neq 1$; or
 \item \label{phi+y}
 $\dim H = 1$,  and every nonzero element of~$\Lie h$ satisfies
$\phi \neq 0$ and $y \neq 0$; or
 \item \label{phi+x}
 $\dim H = 1$,  and every element of~$H$ satisfies $y = 0$ and $\|x\|^2 =
-2 \phi \eta$.
 \end{enumerate}
 \end{enumerate}
 \end{cor}

\subsection{When is the size of $\rho(h)$ quadratic?}
\label{HinN-square-sect}

\begin{lem} \label{HinN-square}
 Assume that $G = \SO(2,n)$, and let $H$ be a closed, connected subgroup of~$N$.
The following are equivalent:
 \begin{enumerate}
 \item \label{HinN-square-seq}
 There is a sequence $\{h_m\}$ of elements of~$H$ with $h_m \to \infty$
and $\rho(h_m) \asymp \|h_m\|^2$.
 \item \label{HinN-square-curve}
 There is a curve $\{h_t\}_{t \in \real^+}$ in~$H$ with $h_t \to \infty$
and $\rho(h_t) \asymp \|h_t\|^2$.
 \item \label{HinN-square-Liealg}
 there is a nonzero element of~$\Lie h$ such that either 
 \begin{enumerate}
 \item \label{HinN-square-indep}
 $\langle (\phi,x), (0,y)\rangle$ is two-dimensional in $\real^{n-1}$,
and either $\dim \Lie h \neq 1$ or $\phi = 0$; or 
 \item \label{HinN-square-x^2=2PhiEta}
 $y = 0$ and $\|x\|^2 = -2 \phi \eta$.
 \end{enumerate}
 \end{enumerate}
 \end{lem}

\begin{proof}
 (\ref{HinN-square-seq} $\Rightarrow$ \ref{HinN-square-Liealg})
 Assume that neither \pref{HinN-square-indep} 
nor~\pref{HinN-square-x^2=2PhiEta} holds. We have $\Lie h
\cap \Lie u _{\alpha +2\beta } = 0$, for otherwise there is an element
of~$\Lie h$ for which $x$, $y$, and~$\phi$ are all zero, so
\pref{HinN-square-x^2=2PhiEta} holds.

\setcounter{case}{0}

\begin{case} \label{case:y=0}
 Assume that $y = 0$ for every element of~$\Lie h$.
 \end{case}
 Because $\Lie h \cap \Lie u_{\alpha+2\beta} = 0$, we know that $\eta$ is
determined by $\phi$ and~$x$, so there exist $p \in \real$ and $b \in
\real^{n-2}$, such that, for every element of~$\Lie h$, we have $\eta =
p\phi + b \cdot x$. 

Because $y = 0$ and $\eta = O \bigl( (\phi,x) \bigr)$, we have $h_{ij} = O
\bigl( (\phi,x) \bigr)$ for all $(i,j) \neq (1,n+2)$. Therefore, if
$\rho(h) \asymp \|h\|^2$, then we must have $h_{1 , n+2} = O \bigl(
(\phi,x) \bigr)$. However, $h_{1 , n+2} = -\phi \bigl( p\phi + b \cdot x
\bigr) - \frac{1}{2} \|x\|^2$ is a quadratic form on~$\Lie h$. Therefore,
because there are vectors $(\phi,x)$ with $h_{1 , n+2}= O \bigl( (\phi,x)
\bigr) \ll \|(\phi,x)\|^2$, this quadratic form must represent~$0$
nontrivially, so there is an element of~$H$ with $h_{1, n+2} = 0$; that
is, $\phi \eta + \|x\|^2/2 = 0$, so
conclusion~\pref{HinN-square-x^2=2PhiEta} holds.

\begin{case} \label{case:phi=0}
 Assume that $\phi = 0$ for every element of~$\Lie h$.
 \end{case}
  From the negation of~\pref{HinN-square-indep}, we see that $x$~and~$y$
are linearly dependent, for every element of~$\Lie h$. Also, we may
assume that there is an element of~$\Lie h$ with $y \neq 0$, for,
otherwise, Case~\ref{case:y=0} applies. Thus, conjugating by an element
of~$U_\alpha$, we may assume that $x = 0$ for every element of~$\Lie h$
\see{conj-alm-same}.  Because $\Lie h \cap \Lie u _{\alpha +2\beta }$,
$\phi$, and~$x$ are all zero, we know that $\eta = O(y)$, so $h_{ij} =
O(y)$ for all $(i,j) \neq (2,n+1)$. On the other hand, $h_{2 , n+1}
\asymp \|y\|^2$. Therefore, $\rho(h) = O(\|y\|^3) \ll \|y\|^4 \asymp
\|h\|^2$.

\begin{case}
 The remaining cases.
 \end{case}
 If $(\phi,x)$ and $(0,y)$ are linearly dependent, for every element
of~$\Lie h$, then either Case~\ref{case:y=0} or Case~\ref{case:phi=0}
applies. Thus, we may assume that there is some element of~$\Lie h$, such
that $\langle (\phi,x) , (0,y) \rangle$ is two-dimensional. (In
particular, we must have $y \neq 0$.) Then, from the negation
of~\pref{HinN-square-indep}, we see that $\dim H = 1$, and $\phi \neq 0$.
Then Lemma~\ref{rho(phiy)} implies that $\rho(h) \asymp \|h\|^{3/2} \ll
\|h\|^2$.

\medbreak

 (\ref{HinN-square-Liealg} $\Rightarrow$ \ref{HinN-square-curve})
 For any element~$h$ of~$U_{\alpha+2\beta}$, we have $h_{ij} = O(1)$ for
all $(i,j) \notin \{ (1,n+1), (2,n+2)\}$, so it is obvious that $\rho(h)
\asymp \|h\|^2$. Thus, we may assume that $\Lie h \cap \Lie u_{\alpha +2\beta }
= 0$.

 Let $h \in H$, and let $h^t = \exp(t \log h)$. 

\setcounter{case}{0}

\begin{case}
 Assume that $\phi = 0$.
 \end{case}
 Because $\phi = 0$, it is clear that $h^t_{ij} = O(t)$ for all $(i,j)
\notin \{ 1,2\} \times \{n+1,n+2\}$. On the other hand, by hypothesis, we
may assume that $\langle (0,x),(0,y)\rangle$ is two-dimensional; that is,
$x$~and~$y$ are linearly independent. Thus,
 $\det
 \left(
 \begin{smallmatrix}
 x \cdot y & \|x\|^2 \\
 \|y\|^2 & x \cdot y \\
 \end{smallmatrix}
 \right)
 \neq 0
 $, so it is easy to see that
  $\det \left(
 \begin{smallmatrix}
 h^t_{1 , n+1} & h^t_{1 , n+2} \\
 h^t_{2 , n+1} & h^t_{2 , n+2} \\
 \end{smallmatrix}
 \right)
  \asymp t^4 \asymp \|h^t\|^2$.
 Therefore, $\rho(h^t) \asymp \|h^t\|^2$.

\begin{case}
 Assume that $y = 0$.
 \end{case}
 In this case, it is clear that 
 $$ h^t_{ij} =
 \begin{cases}
 O(1) & \mbox{if $i \neq 1$ and $j \neq n+2$} \\
 O(t) & \mbox{if $(i,j) \neq (1,n+2)$} \\
 \end{cases}
 .
 $$
 By hypothesis, we may assume that $\|x\|^2 = -2 \phi \eta$, so $h^t_{1
, n+2} = 0$. Thus, because $\phi$, $x$, and~$\eta$ cannot all be zero,
it is clear that $\rho(h^t) \asymp t^2 \asymp \|h^t\|^2$.

\begin{case}
 Assume that both of $\phi$ and $y$ are nonzero. 
 \end{case}
 We see from~\pref{HinN-square-Liealg} that $\dim H \neq 1$, so there is
an element~$u$ of~$H$ with $\phi_u = 0$. We know that $y_u = 0$, for,
otherwise, $[\log u, [\log u, \log h]]$ would be a nonzero element of
$\Lie h \cap \Lie u_{\alpha+2\beta}$. Then, because
$\Lie h \cap \Lie u_{\alpha+2\beta} = 0$, we must have $x_u \neq 0$. 

For each large $t \in \real$, define
 $f_t \colon \real \to \real$ by 
 $$f_t(s) = -t \phi_h ( t \eta_h + s \eta_u)
 - \frac{1}{2} \| t x_h + s x_u \|^2
 + \frac{1}{24} t^4 \phi_h^2 \|y_h\|^2 .$$
 Because $t$ is large, we know that $f_t(0) > 0$, so there is some $s =
s(t) \in \real^+$ with $f_t(s) = 0$. Note that $s$ is approximately
$\frac{1}{\sqrt{12}} t^2 |\phi_h| \|y_h\|/\|x_u\|$.

Let $h^t = \exp( t \log h + s \log u)$. Then $h^t_{1,n+2} = f_t(s) = 0$,
so we see that $h \asymp t^3$. Also, we have
 $$ \det \begin{pmatrix}
 h^t_{1,n+1} & h^t_{1,n+2} \\ 
 h^t_{2,n+1} & h^t_{2,n+2} \\ 
 \end{pmatrix}
 = h^t_{1,n+1} h^t_{2,n+2} \asymp t^6 \asymp \|h^t\|^2 ,$$
 as desired.
 \end{proof}

In the statement of the following corollary, we use
Proposition~\ref{SO1n-whenconj} to convert the condition that $\|x\|^2 \neq
-2 \phi \eta$ to the condition that $\|b\|^2 - \|c\|^2 - 2p < 0$.

\begin{cor}[of proof] \label{HinN-nosquare}
 Assume that $G = \SO(2,n)$, and let $H$ be a closed, connected subgroup of~$N$.
The following are equivalent:
 \begin{enumerate}
 \item
 There is no sequence $\{h_m\}$ of elements of~$H$ with $h_m \to \infty$
and $\rho(h_m) \asymp \|h_m\|^2$.
 \item
 There is no curve $\{h_t\}_{t \in \real^+}$ in~$H$ with $h_t
\to \infty$ and $\rho(h_t) \asymp \|h_t\|^2$.
 \item
 Either
 \begin{enumerate}
 \item \label{dim<x,y>=1}
 for every nonzero element of~$\Lie h$, we have $\phi = 0$ and $\dim
\langle x,y \rangle = 1$; or
 \item \label{b^2+2a<0}
 there exists a subspace~$X_0$ of~$\real^{n-2}$, $b \in X_0$, $c \in
X_0^\perp$, and $p \in \real$ with $\|b\|^2 - \|c\|^2 - 2p < 0$, such that
for every element of~$\Lie h$, we have $y= 0$, $x \in \phi c + X_0$, and
$\eta = p\phi + b \cdot x$; or
 \item 
 $\dim H = 1$, and we have $\phi \neq 0$ and $y \neq 0$
for every element of~$H$; or
 \item \label{c^2+2p<>0}
 $\dim H = 1$ and we have $y = 0$ and $\|x\|^2 \neq -2\phi\eta$
for every element of~$H$.
 \end{enumerate}
 \end{enumerate}
 \end{cor}

\subsection{Some calculations of $\mu(H)$} \label{HinN-calculate}

\begin{lem} \label{SO2n-mu(inv)=mu(h)}
 Assume that $G = \SO(2,n)$. Then $\mu(g^{-1}) = \mu(g)$ for all $g \in G$.
 \end{lem}

\begin{proof} This follows from the fact that there is an element of the
Weyl group of $\SO(2,n)$ that sends each element of~$A$ to its inverse,
so each element of~$A$ is conjugate to its inverse under~$K$.
 \end{proof}

\begin{lem} \label{dim<>1-square}
 Assume that $G = \SO(2,n)$, and let $H$ be a connected, closed subgroup
of~$N$.
 Assume that for every nonzero element of~$\Lie h$, we have $\phi = 0$
and $\dim \langle x,y \rangle \neq 1$. Then $\rho(h) \asymp \|h\|^2$,
 for every $h \in H$.
 \end{lem}

\begin{proof}
 For future reference we prove the stronger fact that for every real
number $t \ge 1$ and every $h \in H$, we have $\rho(ah) \asymp \|ah\|^2$,
where 
 $$a = \diag(t,t,1,1,\ldots,1,1,t^{-1},t^{-1}) .$$
 Let $g = ah$.
 We have
 $ \det 
 \begin{pmatrix}
 g_{1,1} & g_{1,2} \\
 g_{2,1} & g_{2,2} \\
 \end{pmatrix}
 = t^2$
 and
 \begin{eqnarray*}
 4
 \left| \det 
 \begin{pmatrix}
 g_{1,n+1} & g_{1,n+2} \\
 g_{2,n+1} & g_{2,n+2} \\
 \end{pmatrix}
 \right|
 &=&
 4t^2
 \left| \det 
 \begin{pmatrix}
 h_{1,n+1} & h_{1,n+2} \\
 h_{2,n+1} & h_{2,n+2} \\
 \end{pmatrix}
 \right| \\
 &=& 4 t^2 \Bigl( \eta^2 + \bigl( \|x\|^2 \, \|y\|^2 - (x \cdot y)^2 \bigr)
\Bigr)
 .
 \end{eqnarray*}
 Because $x$~and~$y$ are linearly independent for all nontrivial $h \in
H$, we have 
 $$\|x\|^2 \, \|y\|^2 - (x \cdot y)^2 \asymp \|x\|^2 \, \|y\|^2
 \asymp \|x\|^4 ,$$
 so we
conclude that $\rho(ah) \asymp \max(t^2, t^2 \eta^2, t^2 \|x\|^4) \asymp
\|a h\|^2$.
 \end{proof}

\begin{notation}
 We realize $\SO(1,n)$ as the stabilizer in $\SO(2,n)$ of the vector
 $$ (0,1,0,\ldots,0,0,-1,0) \in \real^{n+2} .$$
 \end{notation}

\begin{prop} \label{SO1n-whenconj}
 Assume that $G = \SO(2,n)$. Suppose that $H$ is a closed, connected subgroup of~$N$
and that there exist a subspace~$X_0$ of~$\real^{n-2}$, vectors $b \in
X_0$ and $c \in X_0^\perp$, and a real number~$p$, such that 
 $$ H = \{\, h \in N \mid x \in \phi c + X_0, ~ \eta = p\phi + b \cdot
x, ~ y = 0 \, \} .$$ 
 If $\dim H \ge 2$, then the following are equivalent:
 \begin{enumerate}
 \item \label{SO1n-conj}
 $H$ is conjugate to a subgroup of $\SO(1,n)$.
 \item \label{SO1n-b2+2a<0}
 We have $\|b\|^2 - \|c\|^2 - 2p < 0$.
 \item \label{SO1n-rho(h)=h}
 $\rho(h) \asymp h$ for every $h \in H$.
 \item \label{SO1n-x^2<>-2phieta}
 We have $\|x\|^2 \neq -2\phi\eta$, for every nonzero element of~$\Lie h$.
 \end{enumerate}
 \end{prop}

\begin{proof}
 (\ref{SO1n-b2+2a<0} $\Rightarrow$ \ref{SO1n-conj}) Suppose that
$\|b\|^2 - \|c\|^2 - 2p < 0$. Then $t:=\sqrt{-(\|b\|^2 -\|c\|^2 - 2p)}$
is a positive real number. Set
 $$g =\begin{pmatrix}
 t & 0 & 0 & 0 & 0 \\
 &t^{-1}& t^{-1}(c-b) &-\frac{1}{t} t^{-1}\|c-b\|^2 & 0 \\
  &  & 0 & (c-b)^T &0 \\
 && &t^{-1}&0\\ 
 &&& &t
 \end{pmatrix} .$$
 Then a matrix calculation, using the facts that $c
\cdot b = 0$ and $c \cdot X_0 = 0$, verifies that $y=0$ and $\eta = \phi$
for every element of $gHg^{-1}$, so $gHg^{-1}\subset \SO(1,n)$.

 (\ref{SO1n-conj} $\Rightarrow$ \ref{SO1n-rho(h)=h}) Because $\SO(1,n)$ is
reductive and $A^+ \cap \SO(1,n)$ is the positive Weyl chamber of a maximal
split torus of~$\SO(1,n)$, we have $\mu\bigl(\SO(1,n)\bigr) \subset A^+
\cap \SO(1,n)$. Then, because 
 $$A \cap \SO(1,n) = \{\, \diag(a,1,1,\ldots,1,1,a^{-1}) \mid a >0 \,\} ,$$
 it is clear that $\rho(g) \asymp g$ for every $g \in \SO(1,n)$.

 (\ref{SO1n-rho(h)=h} $\Rightarrow$ \ref{SO1n-x^2<>-2phieta})
 This is immediate from Lemma~\fullref{HinN-square}{-x^2=2PhiEta}.

 (\ref{SO1n-x^2<>-2phieta} $\Rightarrow$ \ref{SO1n-b2+2a<0})
 We prove the contrapositive.

 Assume for the moment that $b \neq 0$. The equation
 $\|sc + b\|^2 = -2s(ps + \|b\|^2)$ is quadratic (or linear) in~$s$.
Because the discriminant $4 \|b\|^2 (\|b\|^2 - \|c\|^2 - 2p)$
of this quadratic equation is nonnegative, the equation has a real
solution $s = \phi_0$. Then, for the element of~$\Lie h$ with
 $\phi = \phi_0$, $x = \phi_0c + b$, $y = 0$, and $\eta = p \phi_0 +
\|b\|^2$,
 we have $\|x\|^2 = -2 \phi \eta$. Because $b \neq 0$, not both of $\phi$
and~$x$ can be~$0$, so we have the desired conclusion.

We may now assume that $b = 0$. Because $\dim H \ge 2$, there is a nonzero
element~$x_0$ of~$X_0$. Because 
 $\|c\|^2 + 2p = -(\|b\|^2 - \|c\|^2 - 2p) > 0$,
 the equation $\|sc + x_0\|^2 = -2s(ps)$ has a solution $s = \phi_0$.
Because $x_0 \neq 0$, we know that $\phi_0 \neq 0$. For the element
of~$\Lie h$ with $\phi = \phi_0$, $x = \phi_0 c + x_0$, and $\eta = p
\phi_0$, we have $\|x\|^2 = -2 \phi \eta$.
 \end{proof}

\begin{lem} \label{rho(phiy)}
 Assume that $G = \SO(2,n)$. Let $H$ be a closed, connected subgroup
of~$N$, such that $\dim H = 1$, and we have $\phi \neq 0$ and $y \neq 0$
for every element of~$H$. Then $\rho(h) \asymp \|h\|^{3/2}$, for every
$h \in H$.
 \end{lem}

\begin{proof}
 Let $h^t = \exp(t u)$, where $u$ is some nonzero element of~$\Lie h$.
 For any large real number~$t$, we see that $h^t_{1,n+2} \asymp t^4$,
but $h^t_{i,j} =O(t^3)$ for $(i,j)
\neq (1,n+2)$. Thus, $h^t \asymp t^4$. Furthermore, we have
 $\det \left( \begin{smallmatrix}
 h^t_{1,n+1} & h^t_{1,n+2} \\
 h^t_{2,n+1} & h^t_{2,n+2} \\
 \end{smallmatrix} \right)
 \asymp t^6$,
 and we have $h^t_{i,j} =O(t^2)$ whenever $i \neq 1$ and $j \neq n+2$,
so we conclude that $\rho(h^t) \asymp t^6 \asymp \|h^t\|^{3/2}$, as
desired.

 For future reference, we note the stronger fact that for every real
number $s \ge 1$ and every $h \in H$, we have $\rho(ah) \asymp
\|ah\|^{3/2}$, where 
 $$a = \diag(s^2,s,1,1,\ldots,1,1,s^{-1},s^{-2}) .$$
 \end{proof}

\section{Subgroups of $\SO(2,n)$ that are not contained in~$N$}
 \label{SO2n-notinN-sect}

\begin{thm} \label{SO2n-semiprod}
 Assume that $G = \SO(2,n)$. Let $H$ be a closed, connected subgroup
of~$AN$, such that $H = (H \cap A) \semiprod (H \cap N)$. The subgroup~$H$
is a Cartan-decomposition subgroup of~$G$ if and only if either
 \begin{enumerate}
 \item $A \subset H$; or
 \item $H \cap N$ is a Cartan-decomposition subgroup of~$G$
\see{SO2n-HinN-CDS}; or
 \item $H \cap A = \ker\alpha$, and we have $\phi=\eta=0$ and $\dim
\langle x,y \rangle = 1$ for every $h \in H \cap N$; or
 \item $H \cap A = \ker \beta$ and we have $y=0$ and $\|x\|^2 = -2
\phi \eta$ for every element of $H \cap N$; or
 \item there is a positive root~$\omega$ such that $\Lie h \cap \Lie n
\subset \Lie u_\omega$, and $H$ satisfies the conditions of
Proposition~\ref{rootsemi-CDS}.
 \end{enumerate}
 \end{thm}

\begin{cor}[of proof] \label{SO2n-semi-notCDS}
 Assume that $G = \SO(2,n)$. Let $H$ be a closed, connected subgroup
of~$AN$, such that $H = (H \cap A) \semiprod (H \cap N)$,
and $H \not\subset N$. The subgroup~$H$ is not a Cartan-decomposition
subgroup of~$G$ if and only if either
 \begin{enumerate}

 \item \label{SO2n-semi-notCDS-1D}
 $H = H \cap A$ is a one-dimensional subgroup of~$A$; or

 \item \label{SO2n-semi-notCDS-<xy>not1}
 $H \cap A = \ker\alpha$, and we have $\phi = 0$ and
$\dim \langle x,y \rangle \neq 1$, for every nonzero element of~$\Lie h
\cap \Lie n$, in which case  $\rho(h) \asymp \|h\|^2$ for all $h \in H$;
or

 \item \label{SO2n-semi-notCDS-y=0}
 $H \cap A = \ker\beta$, and we have $\phi=0$, $y = 0$, and $x \neq
0$, for every nonzero element of $\Lie h \cap \Lie n$, in which case
$\rho(h) \asymp h$ for all $h \in H$; or

 \item \label{SO2n-semi-notCDS-x=0}
 $H \cap A = \ker(\alpha+\beta)$, and we have $\phi=0$, $x = 0$,
and $y \neq 0$, for every nonzero element of $\Lie h \cap \Lie n$, in
which case $\rho(h) \asymp h$ for all $h \in H$; or

 \item \label{SO2n-semi-notCDS-b2-2p}
 $H \cap A = \ker \beta$, and there exist a
subspace~$X_0$ of~$\real^{n-2}$, $b \in X_0$, $c \in X_0^\perp$, and $p
\in \real$, such that $\|b\|^2 - \|c\|^2 - 2p < 0$, and we have $y = 0$,
$x \in \phi c + X_0$, and $\eta = p\phi + b \cdot x$ for every $h \in H$,
in which case $\rho(h) \asymp h$ for every $h \in H$; or

 \item \label{SO2n-semi-notCDS-3/2}
 $H \cap A = \ker(\alpha-\beta)$, $\dim H = 2$,  and there exist
${\hat \phi} \in \Lie u_\alpha$ and ${\hat y} \in \Lie
u_{\beta }$, such that ${\hat \phi} \neq 0$, ${\hat y} \neq 0$, and $\Lie
h \cap \Lie n = \real({\hat \phi} + {\hat y})$, in which case
  $\rho(h) \asymp \|h\|^{3/2}$ for every $h \in H$; or

 \item \label{SO2n-semi-notCDS-dim2}
 $H \cap A =\ker\beta$, $\dim H = 2$, and we have $y = 0$
and $\|x\|^2 \neq -2 \phi \eta$ for every $h \in H$, in which case
$\rho(h) \asymp h$ for every $h \in H$; or

 \item \label{SO2n-semi-notCDS-root}
 there is a positive root~$\omega$ and a one-dimensional
subspace~$\Lie t$ of~$\Lie a$, such that $\Lie h \cap \Lie n \subset \Lie
u_\omega$, $\Lie h = \Lie t + (\Lie h \cap \Lie n)$, and
Proposition~\ref{rootsemi-CDS} implies that $H$ is not a
Cartan-decomposition subgroup.

 \end{enumerate}
 \end{cor}

\begin{proof}[{\bf Proof of Theorem~\ref{SO2n-semiprod}.}] We may assume
that $\Lie h \cap \Lie n$ is not contained in any root space, for
otherwise Proposition~\ref{rootsemi-CDS} applies. We may also assume
that $H \not\subset A$, and that $H \cap N$ is not a
Cartan-decomposition subgroup. The proof proceeds in cases, suggested by
Theorem~\ref{HinN-notCDS}.

\setcounter{case}{0}

\begin{case} \label{semiprod:phi=0,dim<>1}
 Assume that $\phi = 0$ and that $\dim \langle x,y \rangle \neq 1$, for
every nonzero element of~$\Lie h \cap \Lie n$.
 \end{case}
 We show that $\rho(h) \asymp \|h\|^2$, for every $h \in H$, so $H$ is
not a Cartan-decomposition subgroup.
 Let $h \in H$, and write $h = au$ with $a \in H \cap A$ and $u \in H
\cap N$. Because $\Lie h \cap \Lie n \not \subset \Lie u_{\alpha
+2\beta}$ (recall that $\Lie h \cap \Lie n$ is not contained in any root
space), we know that $\Lie h \cap \Lie n$ projects nontrivially into
$\Lie u_\beta + \Lie u_{\alpha +\beta}$. This projection is normalized
by $H \cap A$, but, from the assumption on $\dim \langle x,y \rangle$,
we know that this projection intersects neither $\Lie u_\beta $
nor~$\Lie u_{\alpha +\beta}$, so we must have $\alpha (H \cap A) = 1$.
Thus, $a$ is of the form $a = \diag(t,t,1,\ldots,1,t^{-1} ,t^{-1} )$. We
may assume that $t \ge 1$, by replacing $h$ with~$h^{-1} $ if necessary
\see{SO2n-mu(inv)=mu(h)}. Then the proof of
Proposition~\ref{dim<>1-square} implies that $\rho(h) \asymp \|h\|^2$.

\begin{case}\label{semiprod:phi=0,dim=1}
 Assume that $\phi = 0$ and that $\dim \langle x,y \rangle = 1$, for every
element of~$\Lie h \cap \Lie n$.
 \end{case}
 We show that $H$ is a
Cartan-decomposition subgroup if and only if $\alpha (H \cap A) = 1$ (and
$\Lie h \cap \Lie n  \subset \Lie u_\beta + \Lie u_{\alpha +\beta}$).
Furthermore, when $H$ is not a Cartan-decomposition subgroup, we show that
$\rho(h) \asymp h$ for every $h \in H$.

($\Leftarrow$) Because $\alpha (H \cap A) = 1$, we know that every element
of $H \cap A$ is of the form 
 $$a = \diag(t,t,1,\ldots,1,t^{-1} ,t^{-1} ) ,$$
  so $\rho(a) \asymp t^{\pm 2} \asymp \|a\|^2$.  From
Lemma~\ref{HinN-linear} we also know that there is a path $h^t \to \infty$
in $H \cap N$ with $\rho(h^t) \asymp h^t$. Thus, $H$ is a
Cartan-decomposition subgroup.

($\Rightarrow$) Suppose that $\alpha (H \cap A) \neq 1$. Then, because $\dim
\langle x,y \rangle = 1$ for every $h \in H \cap N$, we see that  $\Lie h
\cap \Lie u_{\alpha +2\beta} = 0$ and either
 \begin{enumerate}
 \renewcommand{\theenumi}{\roman{enumi}} 
 \item $\Lie h \cap \Lie n \subset \Lie u_\beta + \Lie u_{\alpha +2\beta}$;
or
 \item \label{semiprod:phi=0,dim1,alpha+beta}
 $\Lie h \cap \Lie n \subset\Lie u_{\alpha +\beta} + \Lie u_{\alpha
+2\beta}$.
 \end{enumerate}
 Because the Weyl
reflection corresponding to the root~$\alpha$ fixes~$\Lie u_{\alpha
+2\beta}$ and interchanges $\Lie u_\beta$ with~$\Lie u_{\alpha +\beta}$,
we may assume that \pref{semiprod:phi=0,dim1,alpha+beta} holds. Then, by
choosing a large negative value for~$p$, we see that $H$ is a subgroup of
a group of the type considered in Case~\ref{semiprod:b2+2a<0} below, so
$\rho(h) \asymp h$ for every $h \in H$.

\begin{case} \label{semiprod:b2+2a<0}
 Assume that there exist a subspace~$X_0$ of~$\real^{n-2}$, $b \in
X_0$, $c \in X_0^\perp$, and $p \in \real$, with $\|b\|^2 - \|c\|^2 - 2p
< 0$, such that we have $y = 0$, $x \in \phi c + X_0$, and $\eta = p\phi
+ b \cdot x$ for every $h \in H$.
 \end{case}
 We show that $H$ is conjugate to a subgroup of $\SO(1,n)$, so  $\rho(h)
\asymp h$ for every $h \in H$. Therefore, $H$ is not a
Cartan-decomposition subgroup.

  Let $h \in H$, and write $h = au$ with $a \in H \cap A$ and $u \in H \cap
N$. Because $y = 0$ for every element of $\Lie h \cap \Lie n$, we know
that $\Lie h \cap \Lie n \subset \Lie u_\alpha +\Lie u_{\alpha +\beta} +
\Lie u_{\alpha +2\beta}$. Also, because $\|b\|^2 - \|c\|^2 - 2 p < 0$,
we know that $p$ and~$c$ cannot both be~$0$, so $\Lie h \cap \Lie n$ does
not intersect~$\Lie u_\alpha$, and, because $\eta = p\phi + b \cdot x$,
we know that $\Lie h \cap \Lie n$ does not intersect~$\Lie u_{\alpha
+2\beta}$. Thus, because $\Lie h \cap \Lie n \not \subset \Lie u_{\alpha
+\beta}$, and $H \cap A$ normalizes $\Lie h \cap \Lie n$, we see that
$\beta (H \cap A) = 1$, so $H \cap A \subset \SO(1,n)$. Furthermore, $H
\cap A$ centralizes $\Lie u_\beta$, so $H \cap A$ is normalized by the
element~$g$ in the proof of Lemma~\ref{SO1n-whenconj} that
conjugates $H \cap N$ to a subgroup of $\SO(1,n)$. Thus, we conclude that
$g$ conjugates all of~$H$ to a subgroup of $\SO(1,n)$.

\begin{case} \label{semiprod:<2}
 Assume that  $\dim (H \cap N) = 1$,  and that $\phi \neq 0$ and $y
\neq 0$ for every nontrivial element of $H \cap N$.
 \end{case}
 We show that $\rho(h) \asymp \|h\|^{3/2}$ for every $h \in H$, so $H$
is not a Cartan-decomposition subgroup.

 Because $H \cap N$ projects nontrivially into both $\Lie u_\alpha$ and
$\Lie u_\beta$, but does not intersect $\Lie u_\beta $, we must have
$\alpha (a) = \beta (a)$ for every $a \in H \cap A$; that is, $a$ is of the
form $a = \diag(t^2,t,1,\ldots,1,t^{-1} ,t^{-2} )$.
 Because $\dim(H \cap N) = 1$, there exist ${\hat \phi} \in \Lie
u_{\alpha}$,  ${\hat x} \in \Lie
u_{\alpha+\beta}$,  ${\hat y} \in \Lie
u_{\beta}$, and ${\hat \eta} \in \Lie u_{\alpha+2\beta}$, such that
$\Lie h \cap \Lie n = \real( {\hat \phi} + {\hat
x} + {\hat y} + {\hat \eta})$. However, because ${\hat \phi} + {\hat
x} + {\hat y} + {\hat \eta}$ is an eigenvector for every element of
$\ker(\alpha-\beta)$, but the restriction of $\ker(\alpha-\beta)$
to $\Lie u_{\alpha+\beta}$ or~$\Lie u_{\alpha+2\beta}$ is different from
the restriction to each of the other root spaces $\Lie u_\alpha$
and~$\Lie u_\beta$, we see that ${\hat x}$ and~${\hat \eta}$ must both
be zero. Thus, $\Lie h \cap \Lie n = \real( {\hat \phi} + {\hat y})$.
From the proof of Lemma~\ref{rho(phiy)}, we conclude that $\rho(h)
\asymp \|h\|^{3/2}$, as desired.

\begin{case}
 Assume that $\dim (H \cap N) = 1$ and that $\Lie h \cap \Lie n$ is not
contained in any root space. 
 \end{case}
 We show that $H$ is a Cartan-decomposition subgroup if and only if either
 \begin{enumerate}
 \item we have $\phi = 0$, $\eta = 0$, and $\dim \langle x,y \rangle = 1$
for every element of~$\Lie h \cap \Lie n$; or
 \item we have $y=0$ and $\|x\|^2 = -2 \phi \eta$ for every
element of $H \cap N$.
 \end{enumerate}
 We also show that if $H$ is not a Cartan-decomposition subgroup, and
neither Case~\ref{semiprod:phi=0,dim<>1}, Case~\ref{semiprod:phi=0,dim=1},
nor Case~\ref{semiprod:<2} applies to~$H$, then we have $\rho(h) \asymp
h$ for every $h \in H$.

 We may assume that $\phi \neq 0$, for some (and hence every) nonzero element
of $\Lie h \cap \Lie n$, for otherwise  Case~\ref{semiprod:phi=0,dim<>1}
or~\ref{semiprod:phi=0,dim=1} applies. Then we may assume that $y = 0$, for
otherwise Case~\ref{semiprod:<2} applies.

 If $\|x\|^2 = -2 \phi \eta$, then Lemma~\ref{HinN-square} implies that
$\rho(u) \asymp \|u\|^2$ for every $u \in H \cap N$. Furthermore, because
$\beta( H \cap A) = 1$, we know that $\rho(a) \asymp a$ for every $a \in H
\cap A$. We conclude that $H$ is a Cartan-decomposition subgroup.

 On the other hand, if $\|x\|^2 \neq -2 \phi \eta$, it is not difficult to
see that $\rho(h) \asymp h$ for every $h \in H$
(cf.~Case~\ref{HinNcurve-x^2=-2phieta} of  (\ref{HinN-linear-Liealg}
$\Rightarrow$ \ref{HinN-linear-curve}) in the proof of
Lemma~\ref{HinN-linear}).
 \end{proof}

We now consider the case where $H$ is not a semidirect product of the
form $T \semiprod U$, with $T \subset A$ and $U \subset N$.

\begin{thm} \label{SO2n-notsemi-CDS}
 Assume that $G = \SO(2,n)$. Let $H$ be a closed, connected subgroup
of~$AN$ that is compatible with~$A$ \see{compatible}.
 Assume that $H \neq (H \cap A) \semiprod (H \cap N)$ and that $H \cap N$
is not a Cartan-decomposition subgroup.
 The subgroup $H$ is a Cartan-decomposition subgroup if and only if $H$ is
abelian and there exists $\omega \in \{\beta, \alpha+\beta\}$ such
that $A \cap (HN) = \ker \omega$, $\Lie h \cap \Lie u_\omega \neq
0$, and $\Lie h \subset \Lie t + \Lie u_\omega + \Lie
u_{\alpha+2\beta}$, where $\Lie t$ is the Lie algebra of~$A \cap (HN)$.
 \end{thm}

We now describe the Cartan projections of those subgroups that are not
Cartan-decomposition subgroups.

\begin{cor}[of proof] \label{SO2n-notsemi-notCDS}
 
 Assume that $G = \SO(2,n)$. Let $H$ be a closed, connected subgroup of~$G$
that is compatible with~$A$ \see{compatible}, and assume that $H \neq (H
\cap A) \semiprod (H \cap N)$. Then there is a positive root~$\omega$,
and a one-dimensional subspace~$\Lie x$ of $(\ker \omega) + \Lie
u_\omega$, such that $\Lie h = \Lie x + (\Lie h \cap \Lie n)$.

 If $H$ is not a Cartan-decomposition subgroup of~$G$, then either:
 \begin{enumerate}
 \item \label{SO2n-notsemi-notCDS-a,a+b}
  $\omega =\alpha$ and $\Lie h \cap \Lie n \subset \Lie u_{\alpha+\beta}$,
in which case
 $\mu(H) \approx \muH{\|h\|}{\|h\|^2/(\log \|h\|)}$; or
 \item \label{SO2n-notsemi-notCDS-a,a+2b}
  $\omega =\alpha$ and $\Lie h \cap \Lie n \subset \Lie u_{\alpha+2\beta}$,
in which case
 $\mu(H) \approx \muH{\|h\|^2/(\log \|h\|)^2}{\|h\|^2}$; or
 \item \label{SO2n-notsemi-notCDS-a+2b,a}
  $\omega ={\alpha + 2\beta}$ and $\Lie h \cap \Lie n \subset \Lie
u_{\alpha}$, in which case
 $\mu(H) \approx \muH{\|h\|^2/(\log \|h\|)^2}{\|h\|^2}$; or
 \item \label{SO2n-notsemi-notCDS-a+2b,b}
  $\omega ={\alpha + 2\beta}$ and either $\Lie h \cap \Lie n \subset \Lie
u_{\beta}$ or $\Lie h \cap \Lie n \subset \Lie u_{\alpha+\beta}$, in
which case
 $\mu(H) \approx \muH{\|h\|}{\|h\|^2/(\log \|h\|)}$; or
 \item \label{SO2n-notsemi-notCDS-b(a+b),w(a+2b)}
  $\omega \in\{\beta, \alpha + \beta\}$, $\Lie h \cap \Lie n \subset \Lie
u_{\omega} + \Lie u_{\alpha+2\beta}$, and $\Lie h \cap \Lie u_\omega = \Lie
h \cap \Lie u_{\alpha+2\beta} = 0$, in which case
 $\mu(H) \approx \muH{\|h\|}{\|h\|^{3/2}}$; or 
 \item \label{SO2n-notsemi-notCDS-b(a+b),c(a+2b)}
  there is a root~$\gamma$ with $\{\omega, \gamma\} = \{\beta, \alpha +
\beta\}$, $\Lie h \cap \Lie n \subset \Lie u_{\gamma} + \Lie
u_{\alpha+2\beta}$, and $\Lie h \cap \Lie u_{\alpha+2\beta} = 0$, in which
case
 $\mu(H) \approx \muH{\|h\|}{\|h\| (\log\|h\|)^2}$; or 
 \item \label{SO2n-notsemi-notCDS-b(a+b),a+2b}
  $\omega \in \{\beta, \alpha + \beta\}$ and $\Lie h \cap \Lie n = \Lie
u_{\alpha+2\beta}$, in which case
 $\mu(H) \approx \muH{\|h\| (\log\|h\|)}{\|h\|^2}$; or
 \item \label{SO2n-notsemi-notCDS-a+b,a}
   $\omega = \alpha + \beta$ and $\Lie h \cap \Lie n = \Lie u_\alpha$, in
which case
 $\mu(H) \approx \muH{\|h\| (\log\|h\|)}{\|h\|^2}$.
 \end{enumerate}
 \end{cor}

\begin{proof}[{\bf Proof of Theorem~\ref{SO2n-notsemi-CDS}.}] 
 Let $T = A \cap (HN)$. There is a positive root~$\omega$ with $\omega(T)
= 1$, and a one-parameter subgroup~$W$ of~$U_\omega$, such that $\Lie h
\subset \Lie t + \Lie w + (\Lie h \cap \Lie n)$, where $\Lie w$ is the
Lie algebra of~$W$ \see{not-semi}. Note that $T$~and~$W$ are contained in
the Zariski closure of~$H$.

It is not difficult to see that every nonabelian subalgebra of~$\Lie n$
contains $\Lie u_{\alpha+2\beta}$. Therefore, if $\Lie h$ does not
contain $\Lie u_{\alpha+2\beta}$, then $\Lie h \cap \Lie n$ is abelian, and
is centralized by~$W$.

  The proof proceeds in cases. If $\dim (H \cap N) \ge 2$, then
Theorem~\ref{HinN-notCDS} implies that one of Cases
\ref{notsemi-<xy>not1}, \ref{notsemi-<xy>=1}, or~\ref{notsemi-phi<>0}
applies. If $\dim(H \cap N) = 1$, then there are many more possibilities
to consider. Let us see that none of them have been overlooked. Because
of Case~\ref{notsemi-HN=a+2b}, we may assume that $\Lie h \cap \Lie n
\not\subset \Lie u_{\alpha+2\beta}$, so if $\phi = 0$ in $\Lie h \cap
\Lie n$, then either Case~\ref{notsemi-<xy>not1} or
Case~\ref{notsemi-<xy>=1} applies. Thus, we may assume that $\phi \neq
0$ in $\Lie h \cap \Lie n$. Because of Case~\ref{notsemi-HN=a}, we may
assume that $\Lie h \cap \Lie n \neq \Lie u_\alpha$, so $\Lie h \cap
\Lie n$ is not normalized by~$A$. This implies that $\omega = \beta$,
and $\Lie h \cap \Lie n \subset \Lie u_\alpha + \Lie u_{\alpha+\beta} +
\Lie u_{\alpha +2\beta}$, so Case~\ref{notsemi-phi<>0} applies.

\setcounter{case}{0}

\begin{case} \label{notsemi-<xy>not1}
 Assume that $\Lie h \cap \Lie n \not \subset
\Lie u_{\alpha +2\beta}$, and that we have $\phi = 0$ and $\dim \langle
x,y \rangle \neq 1$, for every nonzero element of~$\Lie h \cap \Lie n$.
 \end{case}
 We show that this is impossible.
  Because $\Lie h \cap \Lie n \not \subset
\Lie u_{\alpha +2\beta}$ and $\dim \langle x,y \rangle \neq
1$, there is some $h \in H$ with both $x$~and~$y$ nonzero. Since $T$
normalizes~$H \cap N$, we must have $\alpha(T) = 1$ (that is, $\omega =
\alpha$), so the Zariski closure of~$H$ contains~$U_\alpha$. We have
$[h,U_\alpha] \subset U_{\alpha+\beta} U_{\alpha + 2\beta}$, but, because
$y \neq 0$, we know that $[h,U_\alpha] \not\subset U_{\alpha +
2\beta}$. Thus, $\Lie h \cap \Lie n$ contains an element of the form ${\hat x}
+ {\hat \eta}$, with ${\hat x} \in \Lie u_{\alpha+\beta}$ and ${\hat \eta} \in \Lie
u_{\alpha+2\beta}$, such that ${\hat x} \neq 0$. This contradicts the
assumption of this case.

\begin{case} \label{notsemi-w=a}
 Assume that $\omega = \alpha$.
 \end{case}
 We show that $H$ is not a Cartan-decomposition subgroup.
 Furthermore, $\Lie h \cap \Lie n$ is contained in either $\Lie
u_{\alpha+\beta}$ or~$\Lie u_{\alpha+2\beta}$, and 
 $\mu(H) \approx \muH{\|h\|}{\|h\|^2/(\log\|h\|)^2}$ or 
 $\mu(H) \approx \muH{\|h\|^2/(\log\|h\|)^2}{\|h\|^2}$, respectively.

 Suppose that $\langle \Lie w, \Lie h \cap \Lie n \rangle$ is nonabelian.
We must have $\Lie u_{\alpha + 2\beta} \subset \Lie h$. Then, since
$\langle \Lie w, \Lie u_{\alpha + 2\beta} \rangle$ is abelian, we must
have $\dim (\Lie h \cap \Lie n) \ge 2$.
 \begin{itemize}
 \item If Condition~\ref{HinN-notCDS}\pref{HinN-<xy>not1} holds, then
Case~\ref{notsemi-<xy>not1} applies.
 \item Conditions~\ref{HinN-notCDS}\pref{HinN-<xy>=1}
and~\ref{HinN-notCDS}\pref{HinN-b2-2a} cannot hold, because $\Lie
u_{\alpha+2\beta} \subset \Lie h \cap \Lie n$.
 \end{itemize}
 We conclude that $H \cap N$ is a Cartan-decomposition subgroup, which
contradicts a hypothesis of the theorem.

We now know that $\langle \Lie w, \Lie h \cap \Lie n \rangle$ is
abelian. We must have $\Lie h \cap \Lie n \subset \Lie u_\alpha + \Lie
u_{\alpha+\beta} + \Lie u_{\alpha + 2 \beta}$, the centralizer of~$\Lie
u_\alpha$ in~$\Lie n$. We know that $\ker \alpha$ normalizes $\Lie h \cap
\Lie n$, that the restrictions of $\alpha$, $\alpha+\beta$, and $\alpha +
2\beta$ to $\ker\alpha$ are all distinct, that $\Lie u_\alpha \not\subset
\Lie h \cap \Lie n$, and that $H \cap N$ is not a Cartan-decomposition
subgroup, so we see from Proposition~\ref{HinN-notCDS} that $\Lie h \cap
\Lie n$ must be contained in either $\Lie u_{\alpha + \beta}$ or~$\Lie
u_{\alpha + 2\beta}$. 

\begin{subcase} Assume that $\Lie h \cap
\Lie n \subset \Lie u_{\alpha + \beta}$.
 \end{subcase}
 We have
 $$H = 
 \bigset{
 \begin{pmatrix}
 e^t & t e^t \phi_0 & e^t x & 0  & -e^t \|x\|^2 /2\\
     & e^t     & 0     & 0      & 0 \\
     &         & \Id   & 0      & - x^T \\
     &         &      & e^{-t}      & te^{-t} \phi_0 \\
     &      &      &       & e^{-t} \\
 \end{pmatrix} 
 }{ {t \in \real, \atop x \in X_0}}
 ,$$
 where $X_0$ is a subspace of~$\real^{n-2}$ and $\phi_0$ is a nonzero real
number.

We may assume 
that $t,\|x\| \ge 1$ (see~\ref{SO2n-mu(inv)=mu(h)}
and~\ref{bdd-change}). Then
 $h \asymp \max \{ t e^t, e^t \|x\|^2 \}$ and
 $\rho(h) \asymp e^{2t} \|x\|^2$. 
 Letting $t = 1$ yields $\rho(h) \asymp \|x\|^2 \asymp h$.
 The largest relative value of $\rho(h)$ is obtained by letting $t \asymp
\|x\|^2$, which results in
 $ \rho(h) \asymp t e^{2t} \asymp \|h\|^2/ (\log \|h\|)$.

\begin{subcase} Assume that $\Lie h \cap \Lie n \subset \Lie u_{\alpha +
2\beta}$.
 \end{subcase}
 We have
 $$
 H = 
 \bigset{
 \begin{pmatrix}
 e^t & t e^t \phi_0 & 0 & e^t \eta  & -t e^t \phi_0 \eta \\
     & e^t     & 0     & 0      & -e^t \eta \\
     &         & \Id   & 0      & 0 \\
     &         &      & e^{-t}      & te^{-t} \phi_0 \\
     &      &      &       & e^{-t} \\
 \end{pmatrix} 
 }{t,\eta \in \real}
 ,$$
 where $\phi_0$ is a nonzero real number.

Assuming $t,|\eta| \ge 1$, we have
 $h \asymp t e^t |\eta|$ and
 $\rho(h) \asymp e^{2t} \eta^2$.
 (In calculating $\|\rho(h)\|$, one must note that 
 $\det \begin{pmatrix}
 h_{1,2} & h_{1,n+2} \\
 h_{2,2} & h_{2,n+2} \\
 \end{pmatrix}
 = 0 $.)
 Letting $t = 1$ yields $\rho(h) \asymp \eta^2 \asymp \|h\|^2$.
 The smallest relative value of $\rho(h)$ is obtained by letting $\eta =
1$, which results in
 $ \rho(h) \asymp e^{2t} \asymp \|h\|^2/ (\log \|h\|)^2$.

\begin{case} \label{notsemi-w=a+2b}
 Assume that $\omega = \alpha + 2\beta$.
 \end{case}
 Since the restrictions of $\alpha$, $\beta$, $\alpha+\beta$, and $\alpha
+ 2\beta$ to $\ker(\alpha + 2\beta)$ are all distinct, and $H \cap N$ is
not a Cartan-decomposition subgroup, we see from
Proposition~\ref{HinN-notCDS} that $\Lie h \cap \Lie n$ must be contained
in either $\Lie u_\alpha$, $\Lie u_\beta$, or~$\Lie u_{\alpha + \beta}$.
Because the Weyl reflection $w_\alpha$ interchanges $\Lie u_\beta$
and~$\Lie u_{\alpha + \beta}$, we may assume that $\Lie h \cap \Lie n$ is
contained in either $\Lie u_\alpha$ or~$\Lie u_\beta$. 

\begin{subcase}
 Assume that $\Lie h \cap \Lie n \subset \Lie u_\alpha$.
 \end{subcase}
 We have
 $$H = 
 \bigset{
 \begin{pmatrix}
 e^t & e^t \phi   & 0    & t e^t \eta_0         &  -t e^t \phi \eta_0\\
     & e^{-t}  & 0       & 0         & -t e^{-t} \eta_0 \\
     &         & \Id        & 0            & 0 \\
     &         &          & e^t           & -e^t \phi \\
     &         &          &               & e^{-t} \\
 \end{pmatrix}
 }{t,\phi \in \real},$$
 where $\eta_0$ is a nonzero real number.

Assuming that $t,|\phi| \ge 1$, we have
 $h \asymp t e^t \phi$ and
 $\rho(h) \asymp e^{2t} \phi^2$.
 (Note that 
 $\det \begin{pmatrix}
 h_{1,2} & h_{1,n+2} \\
 h_{2,2} & h_{2,n+2} \\
 \end{pmatrix}
 = 0 $.)
 Letting $t = 1$ yields $\rho(h) \asymp \phi^2 \asymp \|h\|^2$. The
smallest relative value of~$\rho(h)$ is obtained by letting $\phi \asymp
1$, which results in
 $\rho(h) \asymp e^{2t} \asymp \|h\|^2/(\log\|h\|)^2$.

\begin{subcase}
 Assume that $\Lie h \cap \Lie n \subset \Lie u_\beta$.
 \end{subcase}
 We have
 $$H = 
 \bigset{
 \begin{pmatrix}
 e^t & 0       & e^t x    & t e^t \eta_0  &  -e^t \|x\|^2/2\\
     & e^{-t}  & 0        &    0           & -t e^{-t} \eta_0 \\
     &         & \Id        & 0             & -x \\
     &         &          & e^t           & 0 \\
     &         &          &               & e^{-t} \\
 \end{pmatrix}
 }{{t \in \real, \atop x \in X_0}},$$
 where $X_0$ is a subspace of~$\real^{n-2}$, and $\eta_0$ is a nonzero
real number.

Assuming $t,\|x\| \ge 1$, we have
 $h \asymp \max \{ e^t \|x\|^2, t e^t \}$ and
 $\rho(h) \asymp e^{2t} \|x\|^2$.
 Letting $t = 1$ yields $\rho(h) \asymp \|x\|^2 \asymp h$. The
largest relative value of~$\rho(h)$ is obtained by letting $t \asymp
\|x\|^2$, which results in
 $\rho(h) \asymp t e^{2t} \asymp \|h\|^2/(\log\|h\|)$.

\begin{case} \label{notsemi-<xy>=1}
  Assume that we have $\phi = 0$ and $\dim \langle x,y \rangle = 1$, for every
nonzero element of~$\Lie h \cap \Lie n$.
 \end{case}
 Because $\dim \langle x,y \rangle = 1$, we know that $\Lie h$ does not
contain $\Lie u_{\alpha+2\beta}$. Therefore, $\Lie h \cap \Lie n$ is
abelian, and is centralized by the one-parameter subgroup~$W$ of~$\Lie
u_{\omega}$. 

We may assume that $\omega \in \{\beta, \alpha + \beta\}$, for otherwise
Case~\ref{notsemi-w=a} or~\ref{notsemi-w=a+2b} applies.

Note that either $x=0$ for every $h \in H \cap N$, or $y=0$ for every $h
\in H \cap N$. (If there is some $h \in H \cap N$ with both $x$~and~$y$
nonzero, then, because  $\omega \neq \alpha$, both $\Lie h \cap \Lie
u_{\beta}$ and $\Lie h \cap \Lie u_{\alpha+\beta}$ must be nonzero. Then,
since $\dim\langle x,y \rangle = 1$, we see that $\Lie h \cap \Lie n$ is
not abelian, which is a contradiction.)
 Because the Weyl reflection corresponding to the root~$\alpha$ 
fixes~$\Lie u_{\alpha + 2\beta}$ and interchanges $\Lie u_{\beta}$
with~$\Lie u_{\alpha + \beta}$, we may assume that $y = 0$ for every $h \in
H \cap N$. Thus, $\Lie h \cap \Lie n \subset \Lie u_{\alpha + \beta} + \Lie
u_{\alpha + 2\beta}$.

\begin{subcase}
 Assume that $\omega = \alpha + \beta$.
 \end{subcase}
 We have
 $$ H = \bigset{
 \begin{pmatrix}
 1 & 0 & tx_0 + x & b_0\cdot x & - \| tx_0 + x \|^2/2 \\
   & e^t&  0      & 0        & - e^t b_0 \cdot x \\
   &    & \Id     & 0        & - (t x_0 + x)^T \\
   &    &         & e^{-t}   & 0 \\
   &    &         &          & 1 \\
 \end{pmatrix}
 }{ {t \in \real, \atop x \in X_0}}
 ,$$
 where $X_0$ is a subspace of~$\real^{n-2}$, and $b_0$~and~$x_0$ are fixed
vectors in~$\real^{n-2}$ with $x_0 \notin X_0$.

Because $x_0 \notin X_0$, we have $\| tx_0 + x \| \asymp \max\{|t|,
\|x\| \}$. Assuming that $t,\|x\| \ge 1$, we have
 $h \asymp \max\{ e^t, \|x\|^2, e^t |b_0 \cdot x| \}$ and
 $\rho(h) \asymp \max \{ t^2 e^t, e^t \|x\|^2 \}$.
 Letting $t = 1$ yields $\rho(h) \asymp \|x\|^2 \asymp h$.

If there is some nonzero $x_1 \in X_0$ with $b_0 \cdot x_1 = 0$, then
letting $x = x_1$ and $e^t = \|x\|^2$ yields
 $\rho(h) \asymp \|x\|^4 \asymp \|h\|^2$, so $H$ is a Cartan-decomposition
subgroup.

On the other hand, if we have $b_0 \cdot x \neq 0$ for every nonzero $x \in
X_0$, then $b_0 \cdot x \asymp x$ for all $x \in X_0$.  The
largest relative value of~$\rho(h)$ is obtained by letting $e^t \asymp
\|x\|$, which results in
 $\rho(h) \asymp  e^{3t} \asymp \|h\|^{3/2}$.

\begin{subcase}
 Assume that $\omega = \beta$.
 \end{subcase}
 We have
 $$ H = \bigset{
 \begin{pmatrix}
 e^t &0 &   e^t x & e^t b_0\cdot x & -e^t \| x \|^2/2 \\
   & 1  &  t y_0      & -t^2\|y_0\|^2/2        & - b_0 \cdot x \\
   &    & \Id     & -t y_0^T        & - x^T \\
   &    &         & 1   & 0 \\
   &    &         &          & e^{-t} \\
 \end{pmatrix}
 }{ {t \in \real, \atop x \in X_0}}
 ,$$
 where $X_0$ is a subspace of~$\real^{n-2}$, and $b_0$~and~$y_0$ are fixed
vectors in~$\real^{n-2}$ with $y_0 \neq 0$ and $y_0 \perp X_0$.

Assuming $t,\|x\|$ are sufficiently large, we have
 $h \asymp e^t \|x\|^2$ and
 $ \rho(h) \asymp t^2 e^t \|x\|^2$.
 Letting $t = 1$ yields $\rho(h) \asymp \|x\|^2 \asymp h$.
  The largest relative value of~$\rho(h)$ is obtained by letting $\|x\|
\asymp 1$, which results in
 $\rho(h) \asymp  t^2 e^t \asymp \|h\| (\log\|h\|)^2$.

\begin{case} \label{notsemi-phi<>0}
 Assume that $\Lie h \cap \Lie n \subset \Lie u_\alpha + \Lie
u_{\alpha+\beta} + \Lie u_{\alpha+2\beta}$, that $\Lie h \cap \Lie n
\not\subset \Lie u_{\alpha+\beta} + \Lie u_{\alpha+2\beta}$, that $\Lie
u_{\alpha+2\beta} \not\subset \Lie h \cap \Lie n$, and that $\Lie h \cap
\Lie n$ is not normalized by~$A$.
 \end{case}
 Because $\Lie h \cap \Lie n$ is not normalized by~$A$, but is
normalized by~$T$, we must have $\beta(T) = 1$, so $\omega = \beta$. 
Also, because $\Lie u_{\alpha+2\beta} \not\subset \Lie h \cap \Lie n$,
we know that $\Lie h \cap \Lie n$ is centralized by~$\Lie w$. Because the
centralizer of~$\Lie w$ (or any other nontrivial subspace of~$\Lie
u_\beta$) projects trivially into~$\Lie u_\alpha$, this is a
contradiction.

\begin{case} \label{notsemi-HN=a+2b}
 Assume that $\Lie h \cap \Lie n = \Lie u_{\alpha+2\beta}$.
 \end{case}
 We may assume that $\omega \in \{\beta, \alpha + \beta \}$, for otherwise
$\omega=\alpha$, so Case~\ref{notsemi-w=a} applies. Then, because the Weyl
reflection $w_\alpha$ interchanges $\beta$ and~$\alpha+\beta$, but
fixes~$\alpha+2\beta$, we may assume that $\omega = \alpha+\beta$. We have 
 $$H = \bigset{
 \begin{pmatrix}
 1   &    0     & tx_0    & \eta      & -t^2 \|x_0\|^2/2\\
     & e^t     & 0     & 0      & -e^t \eta \\
     &         & \Id     & 0      & -tx_0^T \\
     &         &       & e^{-t} &  0 \\
     &         &       &        & 1 \\
 \end{pmatrix}
 }{ t,\eta \in \real},$$
 where $x_0$ is a nonzero vector in~$\real^{n-2}$. 
 Assuming $t,|\eta| \ge 1$, we have
 $h \asymp e^t |\eta|$ and
 $\rho(h) \asymp \max\{ t^2 e^t, e^t \eta^2 \}$.
 Letting $t = 1$ yields $\rho(h) \asymp \eta^2 \asymp \|h\|^2$.
 The smallest relative value of~$\rho(h)$ is obtained by letting $\eta
\asymp t$, which results in
 $\rho(h) \asymp  t^2 e^t \asymp \|h\| (\log\|h\|)$.

\begin{case} \label{notsemi-HN=a}
 Assume that $\Lie h \cap \Lie n = \Lie u_{\alpha} $.
 \end{case}
 We may assume that $\omega \in \{\beta, \alpha + \beta \}$, for otherwise
$\omega = \alpha + 2\beta$, so Case~\ref{notsemi-w=a+2b} applies. Since
$\Lie u_{\alpha+2\beta} \not\subset \Lie h \cap \Lie n$, we know
that $\Lie w$ centralizes $\Lie h \cap \Lie n$. Then, because $\Lie
u_\alpha$ does not centralize any nontrivial subspace of~$\Lie u_\beta$,
we conclude that $\omega = \alpha + \beta$. We have 
 $$H = \bigset{
 \begin{pmatrix}
 1   &  \phi      & tx_0    &   0    & -t^2 \|x_0\|^2/2\\
     & e^t     & 0     & 0      &0 \\
     &         & \Id     & 0      & -tx_0^T \\
     &         &       & e^{-t} &  -e^{-t} \phi\\
     &         &       &        & 1 \\
 \end{pmatrix}
 }{ t,\phi \in \real},$$
 where $x_0$ is a nonzero vector in~$\real^{n-2}$. 
 Assuming $t,|\phi| \ge 1$, we have
 $h \asymp \max\{e^t, |\phi|\}$ and
 $\rho(h) \asymp \max\{ \phi^2/e^t, t^2 e^t \}$.
 Letting $t = 1$ yields $\rho(h) \asymp \phi^2 \asymp \|h\|^2$.
 The smallest relative value of~$\rho(h)$ is obtained by letting $\phi
\asymp  t e^t$, which results in
 $\rho(h) \asymp  t^2 e^t \asymp \|h\| \log\|h\|$.
 \end{proof}


\begin{thebibliography}{MMM}

\bibitem[Ben]{Benoist}
 Y.~Benoist, Actions propres sur les espaces homog\`enes r\'eductifs,
 \emph{Ann. Math.} 144 (1996) 315--347.

\bibitem[Bor]{Borel-AlgicGrp}
 A.~Borel, \emph{Linear Algebraic Groups, 2nd ed.,} 
 Springer-Verlag, New York, 1991.

\bibitem[BT]{BorelTits}
 A.~Borel and J.~Tits, Groupes r\'eductifs, \emph{Publ. Math. IHES} 27
(1965) 55--150.

\bibitem[Hoc]{Hochschild-Lie}
 G.~Hochschild, \emph{The Structure of Lie Groups},
 Holden-Day, San Francisco, 1965.

\bibitem[Hm1]{Humphreys-LieAlg}
 J.~E.~Humphreys, \emph{Introduction to Lie Algebras and Representation
Theory,}
 Springer-Verlag, New York, 1972.

\bibitem[Hm2]{Humphreys-AlgicGrp}
 J.~E.~Humphreys, \emph{Linear Algebraic Groups,}
 Springer-Verlag, New York, 1975.

\bibitem[Jac]{Jacobson}
 N.~Jacobson, \emph{Lie Algebras}, Dover, New York, 1962.

\bibitem[KPS]{KnightPillaySteinhorn}
 J.~Knight, A.~Pillay and C.~Steinhorn,
 Definable sets in ordered structures~II,
 \emph{Trans. Amer. Math. Soc.} 295 (1986) 593--605.

\bibitem[Kb1]{Kobayashi-isotropy}
 T.~Kobayashi, On discontinuous groups acting on homogeneous spaces with
non-compact isotropy groups,
 \emph{J. Geom. Physics} 12 (1993) 133--144.

\bibitem[Kb2]{Kobayashi-criterion}
 T.~Kobayashi, Criterion of proper actions on homogeneous spaces of
reductive groups,
 \emph{J. Lie Th.} 6 (1996) 147--163.

\bibitem[Kb3]{Kobayashi-survey}
 T.~Kobayashi, Discontinuous groups and Clifford-Klein forms of
pseudo-Riemannian homogeneous manifolds,
 in: B.~\O rsted and H.~Schlichtkrull, eds.,
 \emph{Algebraic and Analytic Methods in Representation Theory},
 Academic Press, New York, 1997, pp.~99--165.

\bibitem[Kos]{Kostant}
 B.~Kostant, On convexity, the Weyl group, and the Iwasawa decomposition,
 \emph{Ann. Sc. ENS.} 6 (1973) 413--455.

\bibitem[Kul]{Kulkarni}
 R.~Kulkarni, Proper actions and pseudo-Riemannian space forms,
 \emph{Adv. Math.} 40 (1981) 10--51.

\bibitem[OW]{OhWitte-CK}
 H.~Oh and D.~Witte,
 Compact Clifford-Klein forms of homogeneous spaces of $\SO(2,n)$
 (preprint).

\bibitem[PS]{PillaySteinhorn}
 A.~Pillay and C.~Steinhorn,
 Definable sets in ordered structures~I,
 \emph{Trans. Amer. Math. Soc.} 295 (1986) 565--592.

\bibitem[Pog]{Poguntke}
 D.~Poguntke, Dense Lie group homomorphisms,
 \emph{J.~Algebra} 169 (1994) 625--647.

\bibitem[Rag]{Raghunathan}
 M.~S.~Raghunathan, \emph{Discrete Subgroups of Lie Groups,}
 Springer-Verlag, New York, 1972.

\bibitem[vdD]{vandenDries}
 L.~P.~D.~van~den~Dries, \emph{Tame Topology and O-minimal Structures},
 Cambridge U. Press, Cambridge, 1998.

\bibitem[Var]{Varadarajan}
 V.~S.~Varadarajan, 
 \emph{Lie Groups, Lie Algebras, and their Representations},
 Springer, New York, 1984. 

\bibitem[W1]{WilkiesJAMS}
 A.~Wilkies, 
 Model completeness results for expansions of the
ordered field of real numbers by restricted Pfaffian
functions and the exponential function,
 \emph{J.~Amer. Math. Soc.} 9 (1996) 1051--1094.

\bibitem[W2]{WilkiesICM}
 A.~Wilkies, 
 O-Minimality,
 \emph{Documenta Mathematica}, Extra Volume ICM~I (1998) 457--460.
 \newline
 {\tt
http://www.mathematik.uni-bielefeld.de/documenta/xvol-icm/99/Wilkie.MAN.html}


\bibitem[Wit]{Witte}
 D.~Witte, Superrigidity of lattices in solvable Lie groups,
 \emph{Invent. Math.} 122 (1995) 147--193.

\bibitem[Zim]{Zimmer}
 R.~J.~Zimmer, \emph{Ergodic Theory and Semisimple Groups},
 Birkh\"auser, Boston, 1984.

\end{thebibliography}
\end{document}